\documentclass[12pt]{amsart}
\usepackage{amscd,amssymb,mathrsfs}

\begin{document}
\title[Crooked surfaces]{Crooked surfaces and anti-de Sitter geometry}

\author[W. Goldman]{William M.\ Goldman}
\address{Department of Mathematics\\
University of Maryland \\
College Park, MD 20742 USA}
\date{\today}
\thanks{We are grateful to the National Science Foundation for
partial support, and in particular the GEAR Network, part of the NSF
program ``Research Networks in the Mathematical Sciences''.}

\newtheorem{thm}{Theorem}
\newtheorem*{MainTheorem*}{Main Theorem}
\newtheorem*{Definition*}{Definition}
\newtheorem{lem}[thm]{Lemma}
\newtheorem{conj}[thm]{Conjecture}
\newtheorem{cor}[thm]{Corollary}
\newtheorem{add}[thm]{Addendum}
\newtheorem{prop}[thm]{Proposition}
\newtheorem{definition}[thm]{Definition}

\newcommand{\R}{\mathbb{R}}
\newcommand{\Z}{\mathbb{Z}}
\newcommand{\C}{\mathbb{C}}
\newcommand{\Norm}{\mathsf{Norm}}
\newcommand{\tr}{\mathsf{tr}}
\newcommand{\Fix}{\mathsf{Fix}}
\newcommand{\Exp}{\mathsf{Exp}}
\renewcommand{\exp}{\mathsf{exp}}
\newcommand{\hExp}{\widehat{\Exp}}
\renewcommand{\det}{\mathsf{det}}
\renewcommand{\dim}{\mathsf{dim}}
\newcommand{\V}{\mathsf{V}}
\newcommand{\E}{\mathsf{E}^3}
\newcommand{\Ad}{\mathsf{Ad}}
\newcommand{\Aut}{\mathsf{Aut}}
\newcommand{\ad}{\mathsf{ad}}
\newcommand{\AdS}{\mathsf{AdS}}
\newcommand{\hAdS}{\widehat{\AdS}}
\newcommand{\SL}{\mathsf{SL}}
\newcommand{\PSL}{\mathsf{PSL}}
\newcommand{\PGL}{\mathsf{PGL}}
\newcommand{\SLtR}{\SL(2,\R)}
\newcommand{\PSLtR}{\PSL(2,\R)}
\newcommand{\PGLtR}{\PGL(2,\R)}
\newcommand{\bs}{\mathbf{s}}
\newcommand{\bS}{\mathbf{S}}
\newcommand{\bJ}{\mathbf{J}}
\newcommand{\sltR}{\frak{sl}(2,\R)}
\renewcommand{\P}{\mathsf{P}} 
\newcommand{\Nn}{\mathfrak{N}} 
\newcommand{\PNn}{\P(\Nn)} 
\newcommand{\Ss}{\mathscr{S}}
\newcommand{\Ein}{\mathsf{Ein}}
\newcommand{\Eto}{\mathsf{E}^{2,1}}
\newcommand{\Rto}{\R^{2,1}}
\newcommand{\bv}{\mathbf{v}}
\newcommand{\bu}{\mathbf{u}}
\newcommand{\bw}{\mathbf{w}}
\newcommand{\SOto}{\mathsf{SO}(2,1)}
\newcommand{\Isom}{\mathsf{Isom}}

\newcommand{\GLtR}{\mathsf{GL}(2,\R)}

\newcommand{\Ht}{\mathsf{H}^2}
\newcommand{\Ho}{\mathsf{H}^1}
\newcommand{\Hn}{\mathsf{H}^n}
\newcommand{\gG}{\mathfrak{g}}
\newcommand{\hH}{\mathfrak{h}}
\newcommand{\kK}{\mathfrak{k}}
\newcommand{\mM}{\mathfrak{m}}
\newcommand{\sS}{\mathfrak{s}}
\newcommand{\bB}{\mathfrak{B}}
\newcommand{\hbB}{\widehat{\mathfrak{B}}}
\newcommand{\bb}{\mathfrak{b}}

\newcommand{\Inn}{\mathsf{Inn}}

\newcommand{\Inv}{\mathsf{Inv}}
\newcommand{\hInv}{\widehat{\mathsf{Inv}}}
\renewcommand{\mod}{\mathsf{mod}}
\newcommand{\Laff}{L^{\mathsf{aff}}}
\newcommand{\Min}{\mathsf{Min}}

\newcommand{\Stem}{\mathsf{Stem}}
\newcommand{\W}{\mathcal{W}} 
\newcommand{\CPM}{\mathcal{C}} 
\newcommand{\hCPa}{\widehat{\mathsf{C}}} 
\newcommand{\CPa}{\mathsf{C}} 
\newcommand{\CS}{\mathfrak{C}} 
\newcommand{\Det}{\mathsf{Det}}

\newcommand{\rpthree}{\R\mathsf{P}^3}
\newcommand{\rpo}{\R\mathsf{P}^1}
\newcommand{\rpn}{\R\mathsf{P}^{n^2-1}}
\newcommand{\Mat}{\mathsf{Mat}}

\newcommand{\one}{\mathbf{1}}
\newcommand{\zero}{{\mathbf{0}}}

\newcommand{\Ker}{\mathsf{Kernel}}
\newcommand{\Image}{\mathsf{Image}}

\newcommand{\hp}{\hat{p}}
\newcommand{\hg}{\hat{g}}
\newcommand{\hh}{\hat{h}}

\newcommand{\CHn}{\mathsf{H}^n_\C}
\newcommand{\CPn}{\mathsf{P}^n_\C}

\newcommand{\bars}{\overline{\sigma}}
\maketitle
\dedicatory{\begin{center}
{\it Dedicated to the memory of Shoshichi Kobayashi}\end{center}}
\begin{abstract}
Crooked planes were defined by Drumm to bound fundamental
polyhedra in Minkowski space for Margulis spacetimes. They were extended by Frances to closed polyhedral surfaces in the conformal compactification
of Minkowski space (Einstein space) which we call crooked surfaces. 
The conformal model of anti-de Sitter space is the interior of the quotient
of Einstein space by an involution fixing an Einstein plane.
The purpose of this note is to show that the crooked planes 
defined in anti-de Sitter space recently by Danciger-Gu\'eritaud-Kassel lift to
restrictions of crooked surfaces in Einstein space which are adapted under
the involution of Einstein space defining anti-de Sitter space.
\end{abstract}

\setcounter{tocdepth}{1} 
\nopagebreak
\tableofcontents
\section*{Introduction}
In 1990, Todd Drumm~\cite{Drumm} 
introduced crooked planes to build fundamental polyhedra for
free discrete groups acting properly and isometrically on $3$-dimen\-sional Minkowski space $\E$. Minkowski space is the Lorentzian analog of 
Euclidean space, and anti-de Sitter space is the Lorenzian analog of hyperbolic
space. 
Recently Jeff Danciger, Fran\c cois Gu\'eritaud and Fanny Kassel~\cite{DGK}
introduced analogous surfaces in $3$-dimensional anti-de Sitter space $\AdS^3$, and have used them to construct
fundamental domains for discrete groups of isometries
of $\AdS^3.$

The conformal compactification of Minkowski space is the 
{\em Einstein universe $\Ein^3$,\/}
the geometry of which has been extensively studied by Charles Frances~\cite{Frances1}.
In particular Frances~\cite{Frances2} considered the extensions of crooked planes to 
$\Ein^3$, which he called {\em generalized crooked planes.\/}
(See 
\cite{BCDGM} 
for an expanded treatment of  Einstein geometry and generalized crooked planes, 
which we renamed {\em crooked surfaces.})
Recently Virginie Charette, Dominik Francoeur, and 
Rosemonde Lareau-Dussault~\cite{CFL}
have used crooked surfaces to build fundamental domains
for Kleinian groups of conformal transformations of 
$\Ein^3$. 
More recently, Danciger-Gu\'eritaud-Kassel have
discussed which complete anti-de Sitter $3$-manifolds
admit fundamental domains bounded by crooked planes.
In particular, this question is more subtle than the case of Margulis spacetimes (in which crooked fundamental
domains always exist). 

The purpose of this note is to interpret the $\AdS$-crooked planes of \cite{DGK} 
in terms of crooked surfaces in $\Ein^3$.
In particular the crooked fundamental polyhedra for 
$\AdS$-manifolds of 
\cite{DGK} will extend to fundamental polyhedra in $\Ein^3$ 
for Lorentzian Kleinian groups as in \cite{CFL}.

%
Here is a precise statement of our main result.
$3$-dimensional Min\-kow\-ski space $\E$ embeds in the $3$-dimensional
Einstein universe $\Ein^3$ as the complement of a lightcone.
The unique nontrivial double covering $\hAdS$ of 
$3$-dimensional anti-de Sitter space $\AdS^3$ embeds in $\Ein^3$
as the complement of a $2$-dimensional Einstein hypersphere 
$\Ein^2\subset\Ein^3$.
The conformal involution $I_\Ss$ 
of $\Ein^3$ whose fixed set equals $\Ein^2$
realizes the deck transformation of the complement 
\[
\Ein^3\setminus \Ein^2  \ = \  \hAdS.
\]
Let $p\in\AdS^3$. 
Let  $\bs\in T_p(\AdS^3)$ be a unit-spacelike vector.
Then the 
{\em crooked plane\/} $\CPa(p,\bs)$ with vertex $p$
and director $\bs$ is defined as the union of all non-spacelike  geodesics tangent to vectors orthogonal to $\bs$ 
and spacelike geodesic {\em rays\/} orthogonal to one of 
the two null directions in $\bs^\perp$, chosen by a procedure
described in \S\ref{sec:Consistent}. 
(See \S\ref{sec:AdSCP} and \cite{DGK} for a detailed definition.) 
A key point, as for classic crooked planes in Minkowski space,
is that a crooked plane in $\AdS^3$ separates $\AdS^3$
into a pair of disjoint {\em crooked halfspaces.}

Lift the basepoint $p\in\AdS^3$ to the double covering
$\hAdS$,  denoting its two preimages under the
double covering $\hAdS\longrightarrow\AdS^3$ by 
\[ \hp, -\hp\ \in\ \hAdS. \] 
Let $\bs\in T_p(\AdS^3)$ be a
unit-spacelike vector, and abusing notation, denote the corresponding
tangent vector in $T_{\hp}(\hAdS)$ by $\bs$ as well.
The {\em tangent cone of $\CPa(p,\bs)$ at $p$\/} is 
the subset of $T_p(\AdS^3)$  consisting of
tangent vectors to smooth paths in $\CPa$ at $p$.
Identify the tangent space $T_p(\AdS^3)$  --- as a Lorentzian vector space ---
with Minkowski space $\E$.

A crooked surface in $\Ein^3$ is determined by an ordered quadruple 
\[
(p_0,p_\infty,p_1,p_2)\in\Ein^3\]
 where $p_0,p_\infty$ are non-incident
and $p_i$ is incident to both $p_0,p_\infty$ for $i=1,2$. 
(Compare \cite{BCDGM,Frances2, LareauDussault,CFL}.)
\begin{Definition*}
The crooked surface $\CS(p_0,p_\infty,p_1,p_2)$ is 
{\em $I_\Ss$-adapted\/} if and only if $I_\Ss$ fixes $p_1$ and $p_2$
and interchanges $p_0$ and $p_\infty$. 
\end{Definition*}
\noindent
This condition is somewhat stronger
than the $I_\Ss$-invariance of  {\em the set\/}    
$\CS(p_0,p_\infty,p_1,p_2)$, 
and rules out potentially interesting hypersufaces which resemble crooked planes, but are not crooked planes. 
\S\ref{sec:AdaptedCS} characterizes this condition in terms of the {\em spine\/} of a crooked surface.

\begin{MainTheorem*}
Let $p\in\AdS^3$ and $\bs\in T_p(\AdS^3)$ be unit-spacelike. 
\begin{enumerate}
\item\label{maintheorem1}
The tangent cone to the $\AdS$-crooked plane 
$\CPa(p,\bs)$ is a crooked plane $\CPM(0_p,\bs)\subset\E\approx T_p(\AdS^3)$ and 
\[
\CPa(p,\bs)\  =\  \Exp_p\big(\CPM(0_p,\bs)\big).
\]
\item\label{maintheorem2}
Let $\hCPa(p,\bs)$ be the inverse image of $\CPa(p,\bs)$ under the double covering 
$\hAdS\longrightarrow\AdS^3$, 
and $\CS(p,\bs)$ be the closure of the image of 
the embedding 
\[
\hCPa(p,\bs) \hookrightarrow \hAdS \xrightarrow{\Psi} \Ein^3.
\]
Then $\CS(p,\bs)$ is an $I_\Ss$-adapted  crooked surface in $\Ein^3$.
\item\label{maintheorem3}
Conversely suppose that $\CS$ is an
$I_\Ss$-adapted  crooked surface in $\Ein^3$.
Then $\CS$ arises from an 
$\AdS$-crooked plane by the above construction.
\end{enumerate}
\end{MainTheorem*}

The first section of the paper develops anti-de Sitter geometry,
defines $\AdS$-crooked planes, describes the double cover, and
proves \eqref{maintheorem1}.
The second section discusses the embedding $\Psi$ 
of the double cover $\hAdS$ in $\Ein^3$ and the conformal model
of anti-de Sitter space.
Using this conformal realization and its symmetries,
the final section proves \eqref{maintheorem2} by computing one example
explcitly. Symmetry is exploited one last time to prove  \eqref{maintheorem3},
the converse statement in the Main Theorem.

Our viewpoint follows the spirit of classical projective geometry, 
but in the {\em conformal\/} context.
Conformal models for hyperbolic and elliptic geometry arise from imposing involutions on the sphere $S^n$, the model space for conformal (Riemannian) geometry.
The conformal model of $\Hn$ consists of the two hemispheres
of $S^n$ complementary to the equator $S^{n-1}\hookrightarrow S^n$, 
each given the Poincar\'e metric. 
Reflection in the equator is the involution $\iota$ {\em defining
hyperbolic geometry,\/} in the following sense. 
Hyperbolic space $\Hn$ is the quotient of the complement of $\Fix(\iota) = S^{n-1}$ by the cyclic group $\langle\iota\rangle$. 
The quotient map
$$
S^n \setminus \Fix(\iota) \longrightarrow
\big(S^n \setminus \Fix(\iota)\big)/\langle\iota\rangle \ =: \ \Hn
$$
is a trivial covering space. Either hemisphere is a fundamental domain
for the action of $\langle\iota\rangle$.

Here anti-de Sitter geometry arises from conformal Lorentzian geometry
(on  the Einstein universe $\Ein^3$) by imposing an involution. Just as the
complement of the fixed-point set of conformal inversion on $S^n$
is a double covering space of hyperbolic space, the complement of
the fixed set of the Lorentzian conformal inversion is a double covering space
of anti-de Sitter space. 
Unlike the Riemannian case above, this covering space is nontrivial.

Metric bisectors in real hyperbolic space $\Hn$ 
are totally geodesic hypersurfaces.
In this model, totally geodesic hypersurfaces are hyperspheres 
orthogonal to the equator. 
Orthogonality to the equator is equivalent to invariance under the reflection
in the equator, the {\em involution defining\/} $\Hn$. 
In this way our observation is a direct analog in the crooked context.

A similar point of view was adopted in \cite{Goldman} to
develop the theory of metric bisectors in complex hyperbolic space
in terms of real analytic hypersurfaces in complex projective space 
$\CPn$ 
which we called {\em extors.\/}
Metric bisectors (or equidistant hypersurfaces) in $\CHn$ analytically contine to 
extors in $\CPn$.
Furthermore they can be characterized as those extors adapted to the
anti-polarity defining $\CHn\subset\CPn$. 
Metric bisectors in {\em complex elliptic space\/} ($\CPn$ with the Fubini-Study
metric) can be similarly characterized as extors in $\CPn$ which are adapted
to the anti-polarity defining complex elliptic geometry.

A consequence of our observation is that the disjointness criterion for
crooked surfaces in $\Ein^3$ proved by Charette, Francoeur, 
and Lareau-Dussault~\cite{CFL}
specializes to a disjointness criterion for the 
$\AdS$-crooked planes defined by 
Danciger, Gu\'eritaud, and Kassel.
(This extends the disjointness criteria for crooked planes developed
in Burelle, Charette, Drumm and Goldman~\cite{BCDG} and 
\cite{DrummGoldman1}).

Since much of the background on this geometry is
not easily accessible in the literature, we have tried
to give a coherent exposition of this geometry.
Thus some parts of the exposition do not directly relate
to the proof of the Main Theorem.
We have tried to indicate these parts, 
including an appendix describing the boundary 
of $\hAdS$.

\subsection*{Acknowledgements}
I thank Jeff Danciger, 
Fran\c cois Gu\'eritaud, 
and Fanny Kassel 
for many helpful conversations, and especially for sharing their 
results~\cite{DGK} with me.
This work clearly owes to Todd Drumm and Charles Frances 
for introducing crooked planes and crooked surfaces.
In addition, I especially thank Virginie Charette,
Dominik Francoeur, and Rosemonde Lareau-Dussault, 
for many illuminating conversations on this subject, 
and in particular for sharing their results \cite{CFL,LareauDussault} with me.

These conversations occurred at 
the special trimester {\em Geometry and Analysis of Surface Group Representations\/} at Institut Henri Poincar\'e in Paris (January-March 2012), 
the GEAR Network Retreat at the University of Illinois in August 2012, the workshop {\em Higher Teichm\"uller-Thurston Theory\/}
at Centre de Recherches Math\'e\-ma\-tiques in Montr\'eal
in October 2012,
and the workshop on {\em Exotic Geometric Structures\/}
and the International Center for Experiment and Research
in Mathematics at Brown University in September 2013.
I am indebted to these venues for their hospitality, 
and in particular for the opportunity to discuss these
mathematical ideas.

I am grateful to Jeff Danciger, Jean-Philippe Burelle and Greg Laun
for reading preliminary versions of this manuscript and suggesting numerous improvements and corrections. 
I thank Thierry Barbot for many illuminating discussions
about Lorentzian geometry.
Finally I would like to thank the anonymous referees
for numerous helpful comments and suggestions.

I dedicate this paper to the memory of Professor Sho
Kobayashi, a good friend who taught me much Differential
Geometry while I was a graduate student at the University of California at Berkeley.

\subsection*{Notation and Terminology}
Homogeneous coordinates on projective space will be denoted
$[X:Y:Z:U:V]$, while inhomogeneous coordinates on affine Minkowski space will
be denoted $(x,y,z)$.
As we predominantly consider Lorentzian $3$-manifolds, with
two spatial dimensions and one temporal dimension,
we shall refer to the {\em ambient\/} dimension:
Thus we denote Minkowski space $\Eto$ by $\E$, 
anti-de Sitter space $\AdS^{2,1}$
by $\AdS^3$ and the Einstein universe $\Ein^{2,1}$ by $\Ein^3$, or simply
$\Ein$. Denote the double cover of $\AdS^3$ by $\hAdS$ and
the $2$-dimensional Einstein universe (one spatial and
one temporal dimension) by $\Ein^2$. 

Denote a symmetric space by $X = G/H$, where $G$ is the isometry
group and $H$ is the isotropy group. 
Denote the exponential map of a complete
pseudo-Riemannian manifold $M$ at a point $p\in M$ by
$$
T_p(M) \xrightarrow{~\Exp_p~} M
$$
and the exponential map for a Lie group $G$ with Lie algebra
$\gG$ by
$$
\gG \xrightarrow{\exp} G.
$$
Denote the adjoint representation of a Lie group on its Lie algebra by
$$
G \xrightarrow{\Ad} \Aut(\gG).
$$
Denote the vector space of $n\times n$ real matrices  by $\Mat_n(\R)$
and the zero and identity matrices by $\zero$ and $\one$, respectively.
Denote the transpose of a matrix $M$ by $M^\dag$.
Denote equivalence classes in $\PSL(2,\R)$ of matrices $A\in\SLtR$ by
$\pm A$. Denote the Lie algebra of $\SLtR$, consisting of traceless
$2\times 2$ real matrices, by $\sltR$, and the hyperbolic plane
(in the Poincar\'e upper halfplane model) by $\Ht$.
Denote the geodesic in $\Ht$ corresponding to the positive imaginary axis by $\Ho$. 
Denote the closure of a subset $S$ by  $\bar{S}$ and 
lifts of points, subsets $S$, etc.\ from a space $M$ to a double covering
space $\widehat{M} \longrightarrow M$ by  $\widehat{S}$.
Denote the set of points fixed by a transformation  $f$ by $\Fix(f)$.
Denote the cyclic group generated by $f$ by $\langle f\rangle$.

If $M$ is a manifold and $p\in M$ is a point, 
denote the zero vector in the tangent space $T_pM$ by
$0_p$.

\section{Anti-de Sitter geometry}
We begin with a review of the geometry of $\AdS$, using the $\PSLtR$-model.
See Danciger~\cite{Danciger}, \S 4.8, for an expanded treatment of this
model and its projective realization.

\subsection{Isometries}
A {\em Lorentzian symmetric space\/} is a Lorentzian manifold $X$ such that
for every point $p$, the isometry $-\one$ of the tangent space $T_pX$ is the
differential of a (necessarily unique) isometry (denoted $\iota_p$) of $X$. 
We call $\iota_p$ the {\em symmetry\/} of $X$ at $p$. 
The isometry group acts transitively. For the general theory of (non-Riemannian)
symmetric spaces, see Eschenburg~\cite{Eschenburg}, 
Helgason~\cite{Helgason}, and Wolf~\cite{Wolf}.

\subsubsection{The $\PSLtR$ model} \label{sec:Isometries}

A convenient model for $3$-dimensional 
anti-de Sitter geometry is the group $\PSLtR$.
Its Lorentzian structure is bi-invariant, 
and  arises from an $\Ad$-invariant inner product
on its Lie algebra $\sltR$. In this case 
$\AdS^3$ corresponds to $\PSLtR$ and the identity component of its isometry
group is 
$$
G^0 := \PSLtR \times\PSLtR
$$ 
acting faithfully by {\em left- and right-multiplication\/} 
on $\PSLtR$:
\begin{equation}\label{eq:LeftRight}
(g_1,g_2):\  
x \longmapsto  g_1 x g_2^{-1}
\end{equation}

The most convenient basepoint in a group is of course its identity
element $e$. 
The isotropy group at $e$ is the subgroup $H\subset G^0$ corresponding to the
diagonal $g_1 = g_2$. Furthermore $H\cong \PSLtR$ acting
by {\em inner automorphisms\/} of $\PSLtR$.

\subsubsection{Components of the isometry group}
Isometries not in the identity component are described, for example,
in Danciger~\cite{Danciger}, \S 4.8 (p.74). 
Briefly, replace
$\PSLtR$ by the group $\PGLtR$ of {\em all\/} 
isometries of $\Ht$, including those which reverse orientation
on $\Ht$.
The identity component of $\PGLtR$ equals $\PSLtR$.
The component consisting of orientation-reversing isometries is
the coset 
$$
R\cdot \PSLtR =  \PSLtR \cdot R \subset \PGLtR
$$
where, for example,
$$
R :=  \pm \bmatrix 1 & 0 \\ 0 & -1\endbmatrix \in \PGLtR
$$
corresponds to a reflection in a geodesic $\Ho:= \Fix(R)\subset\Ht$
(the imaginary axis in the upper halfplane model of $\Ht$).

Now consider the  action of $\PGLtR\times\PGLtR$ by left- and right-multiplication on $\PGLtR$ defined by \eqref{eq:LeftRight}.
Then the group $G$
of all orientation-preserving isometries of $\AdS^3$ equals the index-two
subgroup of  $\PGLtR\times\PGLtR$ which stabilizes the component
$\PSLtR$ of $\PGLtR$.
The orientation-preserving isometry group $G$ contains two components;
the identity component is $G^0 =  \PSLtR\times\PSLtR$.
The new component is the coset 
$$
(R,R)\cdot G^0 = G^0\cdot (R,R) \subset G.
$$ 
This coset consists of isometries which preserve the 
{\em ambient orientation\/} of $\AdS^3$ but reverse {\em time-orientation.\/}

This describes an isomorphism $G\xrightarrow{\cong} \mathsf{PSO}(2,2)$.
Transformations of the form $(Rg_1,g_2)$ or $(g_1,Rg_2)$,
where $g_1,g_2\in\PSLtR$,
do {\em not\/} act orthogonally, but rather take the quadratic form
to its negative. 

\subsubsection{Symmetries in points}

The {\em symmetry\/} of the symmetric space $\AdS^3$ at  $e$ 
is given by group-inversion:
\begin{align}\label{eq:GroupInversion}
\PSLtR &\xrightarrow{\ \iota_e\ } \PSLtR \notag \\
g&\longmapsto g^{-1}.
\end{align}
It is equivariant respecting the permutation automorphism of $
\PSLtR\times\PSLtR$:
\begin{align*}
\PSLtR\times\PSLtR &\longrightarrow \PSLtR\times\PSLtR \\
(g_1,g_2) &\longmapsto (g_2,g_1)
\end{align*}
because \big(applying \eqref{eq:GroupInversion} to \eqref{eq:LeftRight}\big):
\begin{equation}\label{eq:Switch}
\iota_e \circ (g_1,g_2) = (g_2,g_1) \circ \iota_e.
\end{equation}
The differential of $\iota_e$ is the antipodal map $-\one$ on the tangent space $T_e(\AdS^3) \cong \sltR$, and it reverses orientation. 
It generates the other cosets of the full group of isometries of $\AdS^3$, 
including isometries which reverse orientation 
as well as possibly reversing time-orientation.

In this way $\AdS^3$ is a Lorentzian symmetric space with $\iota_e$ as the symmetry in $e$.
More generally, the symmetry in the point $g\in\AdS^3 = \PSLtR$ 
is:
\begin{align}\label{eq:Symmetry}
\AdS^3 &\xrightarrow{\ \iota_g\ } \AdS^3\notag \\
x &\longmapsto  g x^{-1} g
\end{align}
\smallskip
\subsection{Geodesics and totally geodesic surfaces}\label{sec:GeodTotGeodSurf}
Geodesics in $\AdS^3$  are left-cosets of one-parameter subgroups in 
$\PSLtR$.
Since
\[
g\ \exp(t\xi)\  =\    \exp \big(t\  \Ad(g)\xi \big)\ g, 
\]
geodesics correspond to right-cosets of one-parameter
subgroups as well.
$\AdS^3$ is {\em geodesically complete:\/} 
since $\exp(t\xi)$ is defined for all $t\in\R$, 
every geodesic extends indefinitely
in its affine parametrization. Furthermore any two points in $\AdS^3$ are connected by a geodesic.

\subsubsection{Geometry of the tangent space}
Before describing the geodesics and totally geodesic surfaces,
we review basic facts about the Lorentzian vector space $\sltR\cong\Rto$.
Let $\Rto$ denote the $3$-dimensional Lo\-r\-entz\-ian vector space with inner
product associated to the quadratic form:
\begin{align*}
\Rto &\longrightarrow \R \\
\bmatrix x \\ y \\ z \endbmatrix &\longmapsto x^2 + y^2 - z^2
\end{align*}
Under the identification
\begin{align}\label{eq:Identification}
\sltR &\longrightarrow \Rto \notag\\
\bmatrix a & b \\ c & -a \endbmatrix 
&\longmapsto \bmatrix a \\ (b + c)/2 \\ (b - c)/2 \endbmatrix
\end{align}
the inner product in $\Rto$ identifies with the inner product
in $\sltR$ associated to the quadratic form:
\begin{align}\label{eq:QuadraticForm}
\sltR &\longrightarrow  \ \R \notag \\
A &\longmapsto -\det(A) = \frac12 \tr(A^2).
\end{align}
We denote the corresponding (Lorentzian) inner product by $A\cdot B$.
(The Killing form 
\[ A \longmapsto \tr\big(\ad(A)^2\big) \] 
on $\sltR$ equals $-8$ times 
the quadratic form defined by  \eqref{eq:QuadraticForm}.)
One reason for this choice of inner product is that if $A$ denotes one of the
the matrices
$$
\bmatrix 1 & 0 \\ 0 & -1 \endbmatrix,\ 
\bmatrix 0 & 1 \\ 1 & 0 \endbmatrix,\ 
\bmatrix 0 & -1 \\ 1 & 0  \endbmatrix
$$
respectively, then $A\cdot A$ equals $1$, $1$, and $-1$ respectively.

Any inner product on a tangent space to a Lie group extends uniquely to a 
left-invariant pseudo-Riemannian metric. 
An inner product on the tangent
space to the identity element $e$ which is invariant under the {\em adjoint
representation $\Ad$\/} extends to a unique pseudo-Riemannian
metric invariant under {\em both\/}  left- and right-multiplications. 
In this way the inner product corresponding to \eqref{eq:QuadraticForm}
defines a bi-invariant  Lorentzian metric on $\PSLtR$, 
which has constant curvature $-1$. 
This Lorentzian manifold is {\em anti-de Sitter space\/} and denoted $\AdS^3$.

Under the identification \eqref{eq:Identification},
 the {\em Lorentzian cross-product\/}
\begin{align*}
\Rto \times \Rto &\longrightarrow \Rto \\
(A,B) &\longmapsto A\times B
\end{align*}
corresponds to $1/2$ the Lie bracket in $\sltR$.
The {\em scalar triple product\/} 
in $\Rto$ identifies with the $3\times 3$ determinant:
$$
(A \times B) \cdot C\ = \ \mathsf{Det}(A,B,C) 
$$
and satisfies, for example, the vector identity
$$
(A \times B) \cdot C\ = \ (B\times C) \cdot A.
$$
In particular if $A,B\in\sltR$, then $A\times B = \frac12 [A,B]$ 
is orthogonal to both $A$ and $B$:
\begin{equation}\label{eq:tripleproductidentity}
(A \times B) \cdot A\ = \ (A \times B) \cdot B\ = \ 0.
\end{equation}

\subsubsection{Cartan decomposition}\label{sec:CartanDecomposition}
The Cartan decomposition $\gG = \hH \oplus \mM$ of $\AdS^3$ is easy
to describe in the $\PSLtR$-model. Namely, $\gG = \sltR\times\sltR$,
and  \eqref{eq:Switch} implies the Cartan involution is the permutation switching the two summands.
The isotropy algebra $\hH$ is the $1$-eigenspace of the Cartan involution
and is the diagonally embedded subalgebra
\begin{align*}
\sltR &\hookrightarrow \sltR\times\sltR \\
A &\longmapsto (A,A)
\end{align*}
tangent to the isotropy group $H = \Inn\big(\PSLtR\big)$. 
Its orthogonal complement $\mM$  is the $(-1)$-eigenspace of the Cartan
involution and consists of the image of the embedding
\begin{align*}
\sltR &\hookrightarrow \sltR\times\sltR \\
A &\longmapsto (A,-A)
\end{align*}
and models the tangent space 
$$
T_e(\AdS^3) \cong \Rto \cong \sltR.
$$

\subsubsection{Geodesics}
$\AdS^3$ is a {\em Lorentzian symmetric space\/} and we recall
several general facts about symmetric spaces.
 
Choose the identity element $e\in\PSLtR$ as basepoint
in $\AdS^3$ and 
identify its tangent space $T_e\big(\PSLtR\big)$ with the Lie algebra $\sltR$
as above. 

Geodesics through $e$ in $\AdS^3$ correspond to one-parameter subgroups in 
$\PSLtR$. In particular the Lie-theoretic exponential map
$$
\sltR \xrightarrow{\exp} \PSLtR
$$
agrees with the Levi-Civita exponential map:
$$
T_e(\AdS^3) \xrightarrow{\Exp_e} \AdS^3
$$
under the identification of the Lie algebra $\sltR$ with the tangent space:
$$
\sltR \xrightarrow{\cong} T_e\big(\PSLtR\big) 
\xrightarrow{\cong} T_e(\AdS^3).
$$
(Compare \cite{Helgason}, pp.\ 224--226, and \cite{Eschenburg}.
In particular \cite{Helgason}, \S V.6 and Exercise 3 (p.\ 227)
describe symmetric spaces arising from Lie groups with the action
of left- and right-multiplications.)

\subsubsection{Transvections}
Recall in a symmetric space $G/H$, the {\em transvection\/} along a geodesic
$\R\xrightarrow{\gamma}G/H$ 
is the one-parameter group of isometries $\Phi_\gamma(t)\in G$
such that:
\begin{itemize}
\item $\Phi_\gamma(t): \gamma(s) \longmapsto  \gamma(s+t)$;
\item 
Parallel transport along $\gamma$ is the differential
$$
\big(D\Phi_\gamma(t)\big)_{\gamma(s)}:
T_{\gamma(s)}(G/H) \longrightarrow T_{\gamma(s+t)}(G/H).
$$
\end{itemize}
As in \S\ref{sec:CartanDecomposition}, the isotropy algebra 
$\hH$ acting on $\AdS^3$ is the diagonal
$$
\hH \ = \ 
\{ (\xi,\xi) \mid \xi\in \sltR \}$$ and its orthogonal complement $\mM$
equals
$$
\mM \ = \ 
\{ (\xi,-\xi) \mid \xi\in \sltR \}.
$$
Therefore the transvection along the geodesic $\gamma(t)$
corresponding to a one-parameter subgroup $\exp(t\xi)\in G$ equals:
\begin{align}\label{eq:Transvection}
\AdS^3\; &\xrightarrow{\quad\Phi_\gamma(t)\quad} \quad\AdS^3\notag \\
x &\longmapsto\ \exp \bigg(\frac{t}2 \xi\bigg)\  x\  \exp\bigg(\frac{t}2 \xi\bigg)
\end{align}

\subsection{Totally geodesic subspaces and Lie triples}

Identify 
$\mM \cong T_e(\AdS^3)$ with $\sltR$ to realize subspaces of the tangent
space as subspaces of $\sltR$.

Let $\gG = \hH \oplus \mM$ 
be the Cartan decomposition.
Then totally geodesic
subspaces correspond to {\em Lie triples,\/} 
that is, subspaces $\sS\subset \mM$
such that 
\[ [[\sS,\sS],\sS]\subset\sS. \]
Namely, if $S\subset X$ is a totally geodesic subspace
containing the basepoint (corresponding to the coset $eH\in G/H \cong X$),
then $T_{eH}S\subset T_{eH}X \cong \mM$ is a Lie triple.
Conversely, if $\Ss$ is a Lie triple, then 
\[ S = \Exp_{eH}(\Ss) = \exp(\Ss) \]
is a totally geodesic subspace.
(Compare \cite{Eschenburg} and \cite{Helgason},\S IV.7.)

In $\AdS^3$, any $2$-plane $\sS$ in $\mM$ is a Lie triple.
This is easy to see as follows. 
Begin with the following observation:
Since $\dim(\sS)=2$ and $[,]$ is skew-symmetric and bilinear,
$\dim([\sS,\sS])\le 1$.  (In fact $\dim([\sS,\sS])= 1$.) 

There are three types of $2$-planes $\sS\subset\mM\cong\sltR$,
depending on the restriction of the inner product.
We show that each of these is a Lie triple, 
beginning with null $2$-planes.
Null $2$-planes $\sS\subset\sltR$ are  Borel subalgebras, 
each of which is conjugate to the subalgebra 
$\bb\subset\sltR$ of upper-triangular matrices.
Then $[\sS,\sS] \cong [\bb,\bb]$ is a line spanned by a nonzero element (corresponding to a nilpotent upper-triangular matrix in $\sltR$). Furthermore
$$
[[\bb,\bb],\bb] = [\bb,\bb] \subset \bb
$$
as desired.

The remaining two cases follow from the vector identity \eqref{eq:tripleproductidentity}. 
In these two cases, $\sS = v^\perp$ where $v\in\sltR$
is a vector normal to $\sS$. Let $A,B,C\in\Ss= v^\perp$;
we must show that $[[A,B],C]\in\Ss$.
Either $A,B$ are linearly dependent (in which case $[A,B]=0$
and there is nothing to prove)
or $[A,B]$ is nonzero.
In that case, \eqref{eq:tripleproductidentity} implies that $[A,B]$ is 
orthogonal to each $A, B$. Since $A,B$ are linearly independent vectors in
a $3$-dimensional vector space, $[A,B]$ is a nonzero multiple of $v$. 
In either case $\ad[A,B]$ preserves $\Ss$, 
so that $[[A,B],C]\in\Ss$. Explicitly, here are the two cases:
When $\sS$ is spacelike, a normal vector is elliptic, and a 
normal vector $v$ of unit length will satisfy the matrix equation $v^2 = -\one$. 
Furthermore $\ad(v)$ acts on $\Ss = v^\perp$ by rotation of order $4$. 
When $\Ss$ is timelike, a unit length normal vector satisfies $v^2 = \one$
and $\ad(v)$ acts on $\Ss$ by a linear map with eigenvalues $\pm 1$. Thus $\Ss$ is a Lie triple, as desired.

Totally geodesic surfaces through the basepoint
$e\in\AdS^3$ are precisely the images $\exp(\sS)= \Exp_e(\sS)$,
where we identify $T_e(\AdS^3)$ with $\mM$.
Thus,  for any point $g\in\AdS^3$ and every 2-plane $\sS\subset T_g(\AdS^3)$,
\[
S\ =\  \Exp_g(\sS)
\]
is the unique totally geodesic surface $S\ni g$ with
$T_g(S) = \sS$. 
Corresponding to the three types of planes in $\sltR$, 
there are  three types of totally geodesic surfaces of $\AdS^3$.
(Compare \S\ref{sec:embedHt},
\S\ref{sec:NullPlanes}, \S\ref{sec:HtInEin},
\S\ref{sec:GeodTotGeodSurf},  
and \S\ref{sec:TotGeodSurf}.)

\subsubsection{Hyperbolic planes in $\AdS^3$}\label{sec:embedHt}
A point $g\in \AdS^3$ determines a {\em dual\/} 
spacelike totally geodesic surface
$g^*\subset\AdS^3$ which is spacelike, as follows.  

At the origin $e\in\AdS^3$, 
timelike geodesics correspond to elliptic one-parameter subgroups, 
all of the form
$$
\{ h E_\theta h^{-1} \mid \theta\in\R \}
$$
where 
\begin{equation}\label{eq:Elliptic}
E_\theta  := 
\bmatrix \cos(\theta) & -\sin(\theta) \\ \sin(\theta) & \cos(\theta) \endbmatrix.
\end{equation}
An elliptic one-parameter subgroup is {\em periodic\/} in $\SLtR$ with period
$2\pi$ and periodic in $\PSLtR\approx \AdS^3$ with period $\pi$.

For an arbitrary point $g\in\AdS^3$, 
the set of timelike
geodesics from $g$ correspond to 
cosets 
$g hE_\theta h^{-1}$
of elliptic one-parameter subgroups, 
where the parametrization by $\theta$ is chosen so that $E_\theta = e$ 
if and only if
$\theta\in \pi\Z$. 
In particular each timelike geodesic in $\AdS^3$ is closed and has length
$\pi$.

Denote the set of {\em involutions\/} (elements of order two)  
in $\PSLtR$ by $\Inv$.
The inverse image $\hInv$  of $\Inv$ in $\SLtR$ consists of elements
of order four, which have the form
$hE_\theta h^{-1}$ where 
$$
\theta \equiv \pm \pi/2\  (\mod~2\pi).
$$ 
The {\em hyperbolic plane dual to $g$\/} is 
the isometrically embedded $\Ht$ defined as:
\begin{align*}
g^*\,  &:=\,  g\  \Inv\  =\,   \Inv\ g \\ & =\,
\Exp_g\Big( 
\{ \bv\in T_g(\AdS^3) \mid   \bv\cdot\bv = -\pi^2/4 \} \Big).
\end{align*}
Here is an elegant description of this duality.
Let $g\in\AdS^3$.
Then the fixed set $\Fix(\iota_g)$ of the symmetry $\iota_g$ in $g$ 
equals the disjoint union $\{g\} \cup\ g^*$.

Explicitly, $e^*=\Inv$ is the image of the isometric embedding
\begin{align}\label{eq:UpperHalfplane}
\Ht &\hookrightarrow \PSLtR\ \approx\  \AdS^3\notag \\
x + i y &\longmapsto
\pm \frac1{y} \bmatrix   x & -(x^2 +y^2) \\ 1 & -x \endbmatrix
\end{align}
of the upper halfplane $y>0$.
The positive imaginary axis $x = 0 < y$ defines a geodesic
$\Ho\subset\Ht$ where 
\begin{equation}\label{eq:HyperbolicGeodesic}
e^t i \longmapsto {\bJ}(t) := 
\bmatrix 0 & -e^t \\ e^{-t} & 0 \endbmatrix.
\end{equation}

\subsubsection{Indefinite planes in $\AdS^3$}\label{sec:Indefinites}
Totally geodesic indefinite $2$-planes are isometrically embedded copies
of $\AdS^2$, which 
are homeomorphic to M\"obius bands, and therefore nonorientable.

One can see this as follows. 
The isometry group $G^0$ acts transitively on pairs $(g,\sS)$, 
where $g\in\AdS^3$ is a point and $\sS\subset T_g(\AdS^3)$ is an indefinite
$2$-plane. Thus we may assume that $g=e$ and $\sS$ is the subspace
of $T_e(\AdS^3) \cong \sltR$ consisting of matrices of the form
$$
 \bmatrix 0 & * \\ * & 0 \endbmatrix.
$$
In $\SLtR$ this subspace exponentiates to the subset $\mathfrak{S}$
of matrices having
equal diagonal entries with determinant $1$:
\begin{equation}\label{eq:StemPlaneInAdS}
\mathfrak{S}\ :=\  
\Big\{ \bmatrix a & b \\ c & a \endbmatrix \Big| 
a^2 - bc = 1 \Big\}
\end{equation}
which is a one-sheeted hyperboloid. Its image in $\AdS^3$ is the quotient
of this hyperboloid by the antipodal map 
$$
\bmatrix a \\ b \\ c \endbmatrix \longmapsto
\bmatrix -a \\ -b \\ -c \endbmatrix,
$$
with quotient homeomorphic to a M\"obius band. 
Since it is nonorientable, and $\AdS^3$ is orientable, 
it is one-sided in $\AdS^3$.

\subsubsection{Null planes in $\AdS^3$}\label{sec:NullPlanes}
Tangent to degenerate $2$-planes in $T_g(\AdS^3)$ are 
degenerate totally geodesic subspaces, which we
call {\em null planes.\/} A convenient model is the 
{\em Borel subgroup\/}  $\bB = \exp(\bb)$ 
consisting of upper-triangular matrices.
The geodesics through $e$ in this null plane are the one-parameter subgroups:
\begin{equation}\label{eq:UpperTriangularExponential}
\exp
\Bigg( t \bmatrix \alpha & \beta \\ 0 & -\alpha \endbmatrix\Bigg) 
\ =\ \bmatrix e^{t\alpha} & \frac{\beta}{\alpha} \sinh(t\alpha) \\ 0 & e^{-t\alpha} \endbmatrix 
\end{equation}
The commutator subgroup consists of unipotent upper-triangular matrices,
the cosets of which are null geodesics. They define a foliation
tangent to the line field detecting the degeneracy of the 
symmetric $2$-tensor induced by the Lorentzian
structure. The complement of one of these null geodesics is the disjoint
union of two {\em null halfplanes.\/} 
Just as in Minkowski space, the ambient orientation of $\AdS^3$ determines
a way to consistently choose a null halfplane bounded by a null geodesic $\ell$.
(See \S\ref{sec:Consistent}.)        

A null geodesic $\ell$ in $\AdS^3$ lies in a unique null plane $P(\ell)$. 
Given any point $g\in\ell$, the null $2$-plane 
$$
T_g(\ell)^\perp \subset T_g(\AdS^3)
$$
is the unique null plane in the tangent space containing 
the line $T_p(\ell)$. The totally geodesic subspace
$$
P(\ell,g) := \Exp_g\big(T_g(\ell)^\perp\big)
$$ 
is the desired null plane. To see that it is independent
of $g$, choose coordinates so that $P$ consists of the image of 
upper-triangular matrices in $\PSLtR$, and let $\Phi_\ell(t)$ 
denote the transvection along $\ell$, taking $g\in \ell$
to $h\in\ell$. Using the formula \eqref{eq:Transvection}
for transvections, 
$$
\bmatrix 1 & t/2 \\ 0 & 1 \endbmatrix 
\bmatrix a & b \\ 0 & a^{-1} \endbmatrix 
\bmatrix 1 & t/2 \\ 0 & 1 \endbmatrix  \ = \ 
\bmatrix a & b + \Big( \frac{a+a^{-1}}{2}\Big) t \\ 0 & a^{-1} \endbmatrix 
\in P(g,\ell)
$$
so that 
$$
P(h,\ell) = \Phi_\ell(t) P(g,\ell) = P(g,\ell)
$$ as desired.
Henceforth denote $P(g,\ell)$ by just $P(\ell)$.

\subsubsection{Consistent choice of orientations}\label{sec:Consistent}
The complement $P(\ell)\setminus \ell$ decomposes as the disjoint union of two
null halfplanes; here is how to choose one of them, which we call the
{\em wing\/} $\W(\ell)$. 

First assume an orientation on $\AdS^3$, which is given by a volume form
$\Det$ on the tangent space $T_g(\AdS^3)$.
Choose a time-orientation as well (although the result of this
construction turns out to be independent of this choice).
Furthermore choose a future-pointing vector 
$\bv$ tangent to $\bold{\ell}$,
and another future-pointing vector $\bu$ in the complement of $P(\ell)$.
Then, for $p\in\ell$, the {\em wing\/}
\begin{equation}\label{eq:Wing}
\W(\ell) :=\Exp_g\Big(\{\bw\in T_gP(\ell) \mid  \Det(\bv,\bu,\bw) > 0\} \Big)
\end{equation}
is a null halfplane bounded by $\ell$. 
Compare  
\cite{BCDG}, \S 3.2.
\subsection{$\AdS$-crooked planes}\label{sec:AdSCP}
Now we define crooked planes in anti-de Sitter space,
following Danciger-Gu\'eritaud-Kassel~\cite{DGK}. 
We begin by giving several definitions, the first one of which
is analogous to the orginal definition in Minkowski space.
From that we give alternate definitions, in terms of the vertex/spine pair, 
and a dual definition which has no direct analog in 
Minkowski space. 
Then we describe how an $\AdS$-crooked plane $\CPa$
determines a crooked plane $\CPM$ in Minkowski space, 
by a {\em tangent cone construction:\/} 
$\CPM$ is the tangent cone of $\CPa$ at its vertex $g$,
where the tangent space $T_g(\AdS^3)$ is identified with Minkowski space.
This immediately implies that $\CPa = \Exp_g(\CPM)$,
establishing \eqref{maintheorem1} of the Main Theorem. 

The main difference between crooked planes in $\AdS^3$ 
and classic crooked planes 
in Minkowski space is the following:
In $\AdS^3$,  the {\em particles,\/}  
the timelike geodesics  lying on the stem, are compact (homeomorphic to circles), but in Minkowski spce, the particles
are noncompact (lines). 

\subsubsection{Stems, hinges, wings and the spine}
Following \cite{DGK},
given a point $g$ and a unit spacelike vector 
$\bs\in T_g(\AdS^3)$,
we define an {\em $\AdS$-crooked plane,\/} as follows.
The orthogonal complement $\bs^\perp\subset T_g(\AdS^3)$ is
a timelike $2$-plane, and it exponentiates to 
a totally geodesic indefinite plane 
\[ S \ := \ \Exp_g(\bs^\perp) \ \ni \  g.\]
The {\em stem\/} $\Stem(g,\bs)$ is defined as the union of the non-spacelike geodesics in $S$
passing through $g$.  
We call the timelike geodesics lying in $S$ passing through
$g$ the {\em particles,\/}  and the (two) null geodesics
lying on $S$ passing through $g$ the {\em hinges.}

$\Stem(g,\bs)$  is bounded by two null geodesics 
$h_1, h_2$, the hinges. 
Each hinge lies in a unique null plane, as in 
\S\ref{sec:NullPlanes}.
The stem is the union of the two hinges and the particles.

The {\em wings\/} $W_i := \W(h_i)$ are defined as the null halfplanes associated
to  the hinges, as in \eqref{eq:Wing}. The {\em crooked plane\/}
$\CPa(g,\bs)$ is then defined as the union:
\[ \CPa(g,\bs) := \Stem(g,\bs) \cup \W_1 \cup \W_2. \]
The {\em spine\/} is the unique spacelike geodesic contained in
the crooked plane, and equals the image $\sigma$ of:
\begin{align*}
\R &\longrightarrow \AdS^3 \\
t &\longmapsto \Exp_g (t \bs)
\end{align*}
The spine lies on the union of two wings $\W_1 \cup \W_2$,
bisecting each one. 
The vertex $g$ and the spine $\sigma$ completely 
determine the crooked plane.

Since the identity component of the isometry group $G^0$ acts transitively on pairs $(g,\bs)$, 
all $\AdS$-crooked planes are equivalent under the group of orientation-preserving isometries of $\AdS^3$. 

\subsubsection{The dual definition}
The duality between points in $\AdS^3$ and hyperbolic planes,
provides an alternate definition of crooked planes in $\AdS^3$,
which is missing in the classic construction of crooked
planes in Minkowski space.

Dual to the vertex $g$ is a hyperbolic plane 
$g^*\subset\AdS^3$.
The stem $\Stem(g,\bs)$ intersects $g^*$ in a geodesic 
$\ell\subset g^*$.
The pair $(g^*,\ell)$, consisting of a hyperbolic plane 
$g^*$ and a geodesic $\ell\subset g^*$ also suffices to determine the crooked plane.

To recover the vertex-spine pair $(g,\bs)$  
from a pair $(H,\ell)$, where $H$ is a hyperbolic plane containing a geodesic $\ell$, 
take the vertex $g$ to be dual to $H$.  
A unique timelike (indefinite) plane $S$  contains $\ell$ 
and is orthogonal to $H$. 
Since $S$ is orthogonal to $H$, it contains $g$. 
Indeed, $S$ is characterized as the unique spacelike plane containing 
$\ell \cup \{g\}$. 
Thus the data $(g,\pm\bs)$ is equivalent to the data $(H,\ell)$.

\subsubsection{The tangent cone crooked plane}
An $\AdS$-crooked plane with vertex $g$ is
{\em star-shaped\/}  about $g$: that is, it is the 
union of geodesic rays emanating from $g$. 
By definition, its {\em tangent cone\/}  $\CPM\subset T_g(\AdS^3)$ at $g$ consists of all vectors $\bv\in T_g(\AdS^3)$ tangent to smooth rays emanating from $g$. 
(It is a {\em cone\/} in that it is invariant under positive scalar multiplications in the tangent space.) 
In particular,
\[ \CPa\  =\  \Exp_g(\CPM) \]
Moreover,  under the identification of the tangent  space $T_g(\AdS^3)$ with Minkowski space,
$\CPM\subset T_g(\AdS^3)$ identifies with  a crooked plane with vertex $0_g\in T_g(\AdS^3)$. 

\subsection{The double covering space $\hAdS$}
The (unique) double covering space, denoted $\hAdS$, corresponds to 
$\SLtR$, 
consisting of real $2\times 2$ matrices 
$$
\bmatrix a & b \\ c & d \endbmatrix $$
with $ad - bc =1$. 
Since $\hAdS = \SLtR$ is a group, it has a natural basepoint, the identity
matrix $\one\in\SLtR$. 
The deck transformation is the involution given by multiplication by $-\one$,
which we denote by
$$
\hg \longmapsto - \hg.
$$
The action of $G^0$ lifts to the action of the group $\SLtR\times\SLtR$ by
left- and right-multiplication \eqref{eq:LeftRight}, 
whose kernel is the diagonally embedded
central subgroup of order two generated by $(-\one,-\one)$.

In $\hAdS$, the timelike geodesics are all closed and have length $2\pi$.

Unlike $\AdS^3$, not every pair of points in $\hAdS$ is connected by a geodesic.
Let $\one \in\hAdS$ be the identity matrix.
Then the point corresponding to $g\in\SLtR$ is connected by a geodesic to 
$\one$ if and only if $\tr(g) > -2$ or $g$ corresponds to $-\one$. 

The duality between points and hyperbolic planes is somewhat
different as well. Choose a point $\hg\in\hAdS$ covering $g\in\AdS^3$ and
a lift  $\hg^*\subset\hAdS$ of the dual hyperbolic plane $g^*\subset\AdS$.
The preimage in $\hAdS$ of $g^*$ is the disjoint union
$\hg^*\cup (-\hg^*)$, where $\hg^*$ is an (arbitrarily chosen) component.

Then the symmetry $\iota_{\hg}$ in $\hg$ fixes both $\hg$ and $-\hg$
but interchanges the planes $\hg^*$  and $-\hg^*$.
On the other hand, 
\begin{align*}
\hAdS &\longrightarrow \hAdS \\
\hh &\longmapsto  - \iota_{\hg}(\hh)
\end{align*}
is another reflection which fixes each of the two planes $\pm\hg^*$ 
obtained by lifting $g^*$, 
but interchanges the {\em points\/} $\hg$ and $-\hg$.

\subsubsection{Lifting crooked planes}\label{sec:LiftCP}

The inverse image $\hCPa(\hg,\bs)$ of $\CPa(g,\bs)$ 
under the double covering $\hAdS\longrightarrow\AdS^3$ is
connected. Choose one lift of $g$ to call the {\em vertex\/} $\hg$
of $\hCPa(\hg,\bs)$ and call the other one the {\em covertex\/} $\hg' = - \hg$.
The inverse image of each hinge has two components,
one containing the vertex and the other containing the covertex.
Call the component containing the vertex a {\em hinge\/} and
the component containing the covertex a {\em cohinge.\/} 
The inverse image of the stem contains
closed timelike geodesics all of which contain the vertex and the covertex.
Each of the two wings of $\CPa(\hg,\bs)$ has two preimages as well, determined
by which hinge or cohinge bounds them.
We call the lifts of the wings which pass through the vertex {\em wings\/} of the
$\hAdS$-crooked plane, and the lifts of the wings through the covertex the
{\em cowings.\/}
Similarly the spine lifts to two spacelike geodesics; the one
containing the vertex we call the {\em spine\/} of the $\hAdS$ crooked
plane and the one containing the covertex its {\em cospine.\/}

Although $\AdS$-crooked planes are star-shaped about their vertices, their lifts are not. 
If $\CPM$ is the tangent cone of an $\hAdS$ crooked plane
$\hCPa$ at its vertex $\hg$, then  $\hExp_{\hg}(\CPM)$ contains neither  cohinges nor cowings. 
However,
\begin{align*}
\hCPa & \; = \; \hExp_{\hg}(\CPM)\  \cup\   - \hExp_{\hg}(\CPM)  \\
& \; = \;  \hExp_{\hg}(\CPM)\  \cup\   \hExp_{- \hg}(-\CPM) 
\end{align*}
where $-\CPM$ denotes the tangent cone of $\hAdS$ at the covertex
$-\hg$.

\section{Embedding in the Einstein universe}
The double covering space $\hAdS$ of $\AdS^3$ embeds in $\Ein^3$.
This is analogous to the embedding of a double covering space
of $\Ht$ in $S^2$ as the upper and lower hemispheres.  
Reflection in the equator permutes these hemispheres,
acting as the deck involution in this covering space. 
The equator is the fixed set of the reflection, and models conformal
Riemannian geometry in one dimension lower.
Similarly conformal inversion $I_{\Ss}$ in the unit (Lorentzian) sphere is
the deck involution in  $\hAdS$ for the double covering
$\hAdS\to\AdS^3$. The complement of $\hAdS$ in
$\Ein^3$ is exactly the set $\Fix(I_{\Ss})$ of fixed points, which is just
$\Ein^2\subset \Ein^3$, analogous to the Riemannian case.

Minkowski space $\E$ also embeds in $\Ein^3$, 
as an affine patch. Following \cite{BCDGM}, 
we briefly review the geometry of $\Ein$, 
describing it as the conformal compactification of $\E$.
We also describe how to see it topologically as the mapping 
torus of the antipodal map of $S^2$.  
Although some of this exposition is not needed for
the proof of the main theorem, we feel it is useful
in clarifying the properties of the involution $I_{\Ss}$
which is so crucial for this theory. 
See also \cite{Frances1,Frances2,LareauDussault,Rossi}.

The Einstein hypersphere $\Ein_2$, which comprises the
points fixed by $I_{\Ss}$, appears as the boundary of
$\hAdS$. Since we identify $\hAdS$ with $\SLtR$, 
$\partial \hAdS$ corresponds to the boundary of $\SLtR$.
For convenience of the reader, we describe $\partial\SLtR$
in the Appendix, and relate it to $\partial\hAdS$.

\subsection{Einstein space}

Following \cite{BCDGM}, {\em Einstein $n$-space\/} $\Ein^n$ is the projectivized nullcone
$\PNn\subset\P(\R^{n,2})$
in the real inner product space $\R^{n,2}$.  
We follow the notation from \cite{BCDGM} for $n=3$:
points in $\R^{3,2}$ have coordinates $(X,Y,Z,U,V)$ and the 
quadratic form is:
\[Q(X,Y,Z,U,V) := X^2 + Y^2 - Z^2 - UV.\]
Thus the nullcone in $\R^{3,2}$ equals:
\[\Nn := \Big\{(X,Y,Z,U,V)\in\R^{3,2}\   \Big|\  Q(X,Y,Z,U,V) = 0 \Big\}\]
and $\Ein^3$ consists of all $1$-dimensional linear subspaces of $\R^{3,2}$
contained in $\Nn$.
Isotropic $2$-planes in $\R^{3,2}$ projectivize to {\em photons,\/}
which are embedded smooth circles in $\Ein^3$.
The union of all photons through a point $p\in\Ein^3$ is the {\em lightcone\/}
$L(p)\subset\Ein^3$. The complement 
\[ \Min(p) := \Ein^3 \setminus L(p) \]
is the {\em Minkowski patch\/} associated to $p$.

\subsubsection{Minkowski space}
Define, respectively, the {\em improper point\/} and the 
{\em origin \/} in homogeneous coordinates:
\begin{align*}
p_\infty &:= [0:0:0:1:0] \ \in \ \Ein^3\\
p_0  &:= [0:0:0:0:1] \ \in \ \Ein^3.
\end{align*}
The Minkowski patch $\Min(p_\infty)$ is the image of the
embedding of {\em three-dimensional Minkowski space\/} in
$\Ein^3$ (see \cite{BCDGM}):
\begin{align}\label{eq:EmbedMin}
\E &\hookrightarrow \PNn= \Ein^3 \notag\\
(x,y,z) & \mapsto [x: y: z : x^2 + y^2 - z^2 : 1]
\end{align}
In homogeneous coordinates, this Minkowski patch 
$\Min(p_\infty)$ 
is defined by $V\neq 0$, and 
the point $p_0$  corresponds to the {\em origin\/} 
$(0,0,0)\in\E$. 
The ordered pair $(p_\infty,p_0)$ ``makes $\E$ into 
a {\em Lorentzian vector space,\/}'' as follows. 
The subgroup comprising all 
\[ g_{(x,y,z)} :=  \exp \bmatrix
0 & 0 & 0 & 0 & x \\
0 & 0 & 0 & 0 & y \\
0 & 0 & 0 & 0 & z \\
2x & 2y & -2z & 0 & 0 \\
0 & 0 & 0 & 0 & 0 \endbmatrix 
 \]
for $x,y,z\in\R$ acts isometrically on $(\R^{3,2},Q)$
and defines a simply transitive action of the vector group
$\R^3$ on $\E$.
Specifically the orbit $g_{(x,y,z)}(p_0)$ is the point in $\E$ having coordinates $(x,y,z)$. 
This defines the group of translations on 
the affine patch, and $Q$ defines a parallel Lorentzian
metric on $\E$. 

The {\em affine lightcone\/} (the ``usual'' light cone in Minkowski space) 
$\Laff(p_0)$
is the intersection 
\[
\E\cap L(p_0) = L(p_0)\setminus L(p_\infty)
\]
and corresponds to  
\[ \{ (x,y,z)\in\E \mid x^2 + y^2 - z^2 = 0 \}. \]
Its complement in $\E$ equals 
$\E\setminus L(p_0)$, and comprises points
where $x^2 + y^2 - z^2 \neq 0$. 

\subsubsection{Inversion in the unit sphere}
The {\em (Lorentzian) unit sphere\/} $\Ss\subset\E$ 
is the one-sheeted hyperboloid 
\[ \Ss := \{(x,y,z)\in\E \mid x^2 + y^2 - z^2 = 1 \}. \]
Analogous to inversion in Euclidean geometry, 
{\em inversion in $\Ss$\/}:
\begin{align*}
\E\setminus L(p_0)\  &\stackrel{\; I_{\Ss}}\longrightarrow\  
\E\setminus L(p_0) \\
(x,y,z) &\longmapsto \frac1{x^2 + y^2 - z^2} (x,y,z)
\end{align*}
is conformal, and we use this inversion to define a coordinate atlas
on $\Ein^3$. 
%
The involution $I_\Ss$ extends to $\Ein^3$ by mapping $p_0$ to $p_\infty$,
mapping $\Laff(p_0)\setminus\{p_0\}$ 
to $\mathfrak{I}$, and fixing $C_\infty$ pointwise.
In homogeneous coordinates, this involution equals:
 \begin{align}\label{eq:EinInvolution}
\Ein^3 &\xrightarrow{I_{\Ss}} \Ein^3 \notag\\
[X:Y:Z:U:V] &\longmapsto [X:Y:Z: V:U].
\end{align}
The image
\[ \mathfrak{I} \ := \ 
I_{\Ss}\big(\Laff(p_0)\setminus\{p_0\}\big) \ =\ 
L(p_\infty) \cap \Min(p_0) \setminus \{p_\infty\}\]
is called the {\em inverse cone.\/} 
The {\em ideal circle\/} is the intersection
\begin{equation}\label{eq:IdealCircle}
C_\infty := L(p_0) \cap L(p_\infty),
\end{equation}
and consists of points with homogeneous coordinates
$\bmatrix X & Y & Z & 0 & 0 \endbmatrix$
with $X^2 + Y^2 - Z^2 = 0$. 
It (together with the improper point $p_\infty$)
compactifies both the unit-sphere $\Ss$ and the lightcone $\Laff(p_0)$.
Alternatively, $C_\infty$ comprises ideal points of photons through the improper point $p_\infty$).

These four sets define the stratification of $\Ein^3$:
\begin{itemize}
\item (Dimension $3$) The Minkowski patch $\E$ is a $3$-cell;
\item (Dimension $2$) The inverse cone $\mathfrak{I}$
consists of two disjoint annuli;
\item (Dimension $1$) The ideal circle $C_\infty$ is a $1$-sphere;
\item (Dimension $0$) The improper point $p_\infty$ is a $0$-cell.
\end{itemize}
For more details see \cite{BCDG,LareauDussault}.

\subsubsection{A remark on the topology of $\Ein^3$.}
Topologically, $\Ein^3$ is a mapping torus of the antipodal map of $S^2$, as can be seen by writing the quadratic form as 
a difference of two sums of squares:
\[
Q(X,Y,Z,U,V)\  =\  
\Big(X^2 + Y^2 + \bigg(\frac{U-V}2\bigg)^2\Big)
\ -\ 
\Big(Z^2 + \bigg(\frac{U+V}2\bigg)^2\Big).
\]
Then $\Ein^3$ corresponds the quotient of the locus
of
\[
X^2 + Y^2 + \bigg(\frac{U-V}2\bigg)^2 = 
Z^2 + \bigg(\frac{U+V}2\bigg)^2 = 1 \]
(homeomorphic to a product $S^2 \times S^1$
by the antipodal map, which is the above mapping torus.

Because of \eqref{eq:IdealCircle}, the union 
$\Min(p_0)\cup\Min(p_\infty)$ 
equals the complement $\Ein^3\setminus C_\infty$. 
The topological complexity of $\Ein^3$ (it is not a sphere)
implies that one needs at least three Minkowski patches to cover it. 
If $p\in\Ein^3$, then
\[ \Min(p_0)\cup\Min(p_\infty) \cup\Min(p)\ =\ \Ein^3 \]
if and only if $p$ lies on a timelike line through $p_0$.
This is the condition that the {\em Maslov index\/} of the triple
$(p_0,p_\infty,p)$ equals $\pm 1$.


\subsection{Embedding $\hAdS$} The mapping
\begin{align}\label{eq:EmbedSL2R}
\SLtR &\stackrel{\Psi}\hookrightarrow  \PNn \notag\\
\bmatrix a & b \\ c & d \endbmatrix &\longmapsto
[ a-d : b+c : b - c: a + d -2 : a+d +2 ]
\end{align}
embeds $\SLtR$ in $\Ein^3$. 
The identity matrix $\one\in \SLtR$ maps
to the point $p_0$ corresponding to the origin $(0,0,0)$ in Minkowski
space $\E$, under  \eqref{eq:EmbedMin}.
The other central element $-\one\in\SLtR$ maps to $p_\infty$.

\subsubsection{Equivariance under the deck involution}
The double covering $\hAdS \longrightarrow \AdS^3$ corresponds
to the double covering $\SLtR \longrightarrow \PSLtR$.
In this model, the deck involution is multiplication by  $-\one$ in 
$\SLtR$, and $\Psi$ relates this involution to the inversion $I_{\Ss}$
defined by \eqref{eq:EinInvolution}.  To  see this, write:
\begin{align*}
M &:= \bmatrix a & b \\ c & d \endbmatrix, \\ 
\Psi(M) & = [X:Y:Z: U: V] \\
& = [ a-d : b+c : b - c: a + d -2 : a+d +2 ],
\end{align*}
so  \eqref{eq:EmbedSL2R} implies:
\begin{align*}
\Psi(-M) & = [-X:-Y:-Z: -(a+d) -2 : -(a+d) +2 ]\\ & =
[X:Y:Z: (a+d) +2 : (a+d) -2 ] \\ &= [X:Y:Z:V:U].
\end{align*}

%

\subsubsection{The boundary of $\hAdS$}
The double cover $\hAdS$ of $\AdS^3$ is bounded by $\Ein^2$,
which is a $1$-sided surface in $\Ein^3$. 
(Although $\Ein^3$ is nonorientable, $\Ein^2$ is
oreintable.) 
This surface consists of all points in $\Ein^3$
fixed under the involution $I_\Ss$ and the quotient mapping
\[
\Ein^3 \longrightarrow \Ein^3/I_\Ss \ =\  \AdS^3 \cup \Ein^2
\]
is a double covering  branched over $\Ein^2$.
Using the usual compactification of $\AdS^3$ as a closed solid torus
(the $3$-manifold-with-boundary $D^2\times S^1$), $\Ein^3$
is homeomorphic to the identification space of 
$D^2\times S^1$ by an orientation-reversing involution indentifying the boundary to itself.

The standard construction of the boundary $\partial\AdS$
and its double covering $\partial\hAdS$ are reviewed
in the appendix. 


The embedding of $\hAdS\xrightarrow{\Psi}\Ein^3$ and the involution
$I_{\Ss}$ intimately relate. 
The complement of the image of 
$\Psi$ is an Einstein hypersphere $\Ein^2$, 
consisting of fixed points of the involution $I_{\Ss}$.
Clearly points in the image of $\Psi$ are not fixed by $I_{\Ss}$; thus
their homogeneous coordinates satisfy $U \neq V$. Conversely,
\begin{align*}
a & = \frac{2 X + U + V}{V - U} \\
b & = \frac{2(Y + Z)}{V - U} \\
c & = \frac{2(Y - Z)}{V - U} \\
d & = \frac{-2 X + U + V}{V - U} 
\end{align*}
defines an element 
$$
\bmatrix a & b \\ c & d \endbmatrix \in\SLtR
$$ 
mapping to $[X:Y:Z:U:V]$ if $U\neq V$.

The hyperplane in $\R^{3,2}$  defined by $U = V$ has signature $(2,2)$
and therefore its image is a $2$-dimensional Einstein subspace, which we henceforth call $\Ein^2$.
It equals the union 
\[ \Ss\cup C_\infty\ = \ \Fix(I_\Ss). \]
Although the unit sphere $\Ss$ (that is, the one-sheeted
hyperboloid $x^2 + y^2 -z^2 = 1$) disconnects $\E$, its closure
$\Ein^2$ does {\em not\/} disconnect $\Ein^3$.
%
\subsubsection{Ideal points for the embedding of $\hAdS\hookrightarrow\Ein^3$}
The lightcone 
$$
L(p_\infty) = I_{\Ss}\big(L(p_0)\big)
$$ 
is the complement of
Minkowski space $\E$ in $\Ein^3$. 
It is the union of all null geodesics through the improper point 
$p_\infty = I_{\Ss}(p_0)$. 
Since $p_\infty$ corresponds to the matrix $-\one\in\SLtR$, 
the lightcone $L(p_\infty)$
corresponds to the cosets of unipotent one-parameter subgroups through 
$-\one$.
These are precisely the elements of $\SLtR$ having trace $-2$. 

As in \cite{BCDGM}, this lightcone meets $L(p_0)$ in the ideal circle
$C_\infty$. Thus the boundary of $\hAdS$  is the Lorentzian unit sphere
$\Ss\subset \E$ compactified with the ideal circle $C_\infty$. The union
$\Ss\cup C_\infty$ is $\Ein^2$, which is topologically a torus, which is a 
$1$-sided surface in the nonorientable $3$-manifold $\Ein^3$.

\subsubsection{Rulings of $\Ein^2$}\label{sec:Rulings}
The fibers of the two projections
\[ \partial\hAdS \approx S^1 \times S^1 \longrightarrow S^1
\]
correspond to the {\em rulings\/} of $\partial\hAdS = \Ein^2$.
(Compare Danciger~\cite{Danciger}.)
In terms of the affine piece of $\Ein^2$, the Lorentzian unit-sphere
$\Ss$ defined by
$x^2 + y^2 - z^2 = 1$,
these are just the two families of null lines 
\begin{equation}\label{eq:Rulings}
\Ss^+_{\theta}\ :=\ \bigg\{
\bmatrix \cos(\theta) -\eta \sin(\theta) \\ 
\sin(\theta) +\eta \cos(\theta) \\  \eta  \endbmatrix
\bigg| \eta\in\R \bigg\}
\end{equation}
and 
$$
\Ss^-_{\theta}\ :=\ 
\bigg\{
\bmatrix \cos(\theta) +\eta \sin(\theta) \\ 
\sin(\theta) -\eta \cos(\theta) \\  \eta  \endbmatrix
\bigg| \eta\in\R \bigg\}
$$
respectively, where $\theta\in\R/2\pi\Z$ is fixed.
These null lines have  respective ideal points
$$
\partial L^{\pm}_\theta\ := \ 
[0:-\sin(\theta): \cos(\theta):\pm1: 0:0]\  \in\  C_\infty.
$$


\subsection{ $\hAdS$ totally geodesic subspaces in $\Ein^3$}
\label{sec:HtInEin}
To relate  $\AdS$-crooked planes,
crooked surfaces in $\Ein^3$ and classic crooked planes in Minkowski space
$\E$,
we first describe how $\hAdS$ totally geodesic subspaces extend to 
submanifolds of $\Ein^3$.

\subsubsection{Geodesics}\label{sec:Geodesics}
As timelike geodesics in $\hAdS$ correspond to elliptic one-parameter subgroups in $\SLtR$, they are all closed and do not approach the boundary
of $\hAdS$ in $\Ein^3$. 

Null geodesics in $\hAdS$ have a single endpoint in $\partial\hAdS = \Ein^2$.
In our conformal Minkowski model of Einstein space, there are two types,
depending on whether or not the null geodesic $\nu$ passes through the improper point. 
If $p_\infty\in\nu$, then $\nu$ lies on 
$\mathfrak{I}\cup\{p_\infty\}$;
its single endpoint lies on the ideal circle $C_\infty$. 
If $\nu$  does not contain $p_\infty$, then it misses the ideal circle and contains a point
of  $\mathfrak{I}$. Its endpoint lies on $\Ss$. 

The embedding $\Psi$ defined in \eqref{eq:EmbedSL2R} maps an  elliptic one-parameter subgroup 
to a timelike geodesic:
\begin{equation}\label{eq:elliptics}
\bmatrix \cos(\theta) & -\sin(\theta) \\
\sin(\theta) & \cos(\theta) \endbmatrix 
\stackrel{\Psi}\longmapsto 
[0:0: \sin(\theta): 1-\cos(\theta) : 1 + \cos(\theta) ]
\end{equation}
which under \eqref{eq:EmbedMin} corresponds to the $z$-axis $x=y=0$,
with
$$
z = \frac{\sin(\theta)}{1+\cos(\theta)} = \tan(\theta/2),
$$
mapping $-\pi < \theta < \pi$ to $-\infty < z < \infty$.

The embedding $\Psi$ defined in 
\eqref{eq:EmbedSL2R}
maps the (cosets of) unipotent one-parameter subgroups to photons:
\begin{align}\label{eq:unipotents}
\bmatrix 1 & t \\ 0 & 1 \endbmatrix\  &\stackrel{\Psi}\longmapsto 
H_1^+(t) :=  [0:t:t:0:4] \longleftrightarrow (0, t/4, t/4) \notag \\
-\bmatrix 1 & t \\ 0 & 1 \endbmatrix\  &\stackrel{\Psi}\longmapsto 
H_1^-(t) :=  [0:t:t:4:0] \notag \\
\bmatrix 1 & 0 \\ t & 1 \endbmatrix\  &\stackrel{\Psi}\longmapsto 
H_2^+(t) :=  [0:t:-t:0:4] \longleftrightarrow (0, t/4, -t/4)\notag \\
-\bmatrix 1 & 0 \\ t & 1 \endbmatrix\  &\stackrel{\Psi}\longmapsto 
H_2^-(t) :=  [0:t:-t:4:0] 
\end{align}
Their images $h^+_1$,  $h^+_2$,  
correspond to null lines  $x = y -z = 0$ and $x = y -z = 0$
respectively. Inverting in $\Ss$ yield 
$h^-_1$,  $h^-_2$ respectively, which lie on the inverse cone.

Finally, $\Psi$ maps cosets of hyperbolic one-parameter subgroups
in $\SLtR = \hAdS$ 
to spacelike arcs. 
An illustrative model is the intersection of the $x$-axis with the one-sheeted
hyperboloid $\Ss$.
The geodesic in $\hAdS$ emanating from $\one$, which corresponds
to the origin $(0,0,0)\in\E$, is just the line segment 
\begin{equation}\label{eq:UnitInterval}
\{ -1 < x < 1\} \times\{(0,0)\}
\end{equation}
The geodesic in $\hAdS$ emanating from $-\one$, which corresponds to the improper point $p_\infty$, divides into the union of two rays 
\begin{equation}\label{eq:TwoRays}
\{ x <  -1 \}\times\{(0,0)\} \ \cup \ 
\{ 1 < x \}\times\{(0,0)\}, 
\end{equation}
separated by the improper point $p_\infty$ at $\infty$.
These two rays correspond to the components of the complement of $-\one$
in the this coset.
Explicitly:
\begin{equation}\label{eq:Hyperbolics}
\bmatrix e^t  & 0 \\
0 & e^{-t} \endbmatrix 
\stackrel{\Psi}\longmapsto 
[\sinh(t):0: 0: \cosh(t) - 1 : \cosh(t) + 1]
\end{equation}
which under \eqref{eq:EmbedMin} corresponds to the $x$-axis $y=z = 0$,
with
$$
x \ =\  
\frac{\sinh(t)}{\cosh(t) + 1} \ = \ 
\tanh(t/2)\  
 $$
mapping $-\infty< t < \infty$ to $-1 < x < 1$.
The other geodesic corresponds to its $I_\Ss$-image:
\begin{equation}\label{eq:MinusHyperbolics}
- \bmatrix e^t  & 0 \\
0 & e^{-t} \endbmatrix 
\stackrel{\Psi}\longmapsto 
[\sinh(t):0: 0: \cosh(t) + 1 : \cosh(t) - 1]
\end{equation}
which under \eqref{eq:EmbedMin} corresponds to the $x$-axis $y=z = 0$,
with
$$
x \ =\   \frac{\sinh(t)}{\cosh(t) - 1} \ = \ \coth(t/2)
$$
mapping $-\infty< t  < 0 < \infty$ to $x < -1$ and
$ 0 < t < \infty$ to $x > 1$.

%

These two $\hAdS$ geodesics, together with their two endpoints on $\Ss$ constitute a {\em spacelike circle\/} in $\Ein^3$, as defined in \cite{BCDGM}.

\subsubsection{Totally geodesic surfaces}\label{sec:TotGeodSurf}
Now we see how totally geodesic surfaces of $\hAdS$ embed in $\Ein^3$.
Recall three types of totally geodesic surfaces in $\AdS^3$ exist:
definite planes, isometric to $\Ht$; 
degenerate planes, corresponding to cosets of Borel subgroups of $\SLtR$;
indefinite planes, isometric to $\AdS^2$. 

A definite plane is dual to a point $g\in\AdS^3$
as discussed in
\S\ref{sec:embedHt}. Here is a simple example.
The plane $e^*$  dual to the identity element $e\in\PSLtR$
is the set of involutions $\Inv\subset\PSLtR$; see 
\eqref{eq:UpperHalfplane} for an explicit embeddiing of $\Ht$ in $\AdS^3$. 
In $\Ein^3$ the two lifts of this hyperbolic plane are the two components
of the hyperboloid $x^2 + y^2 - z^2 = -1$; their common boundary is the
ideal circle $C_\infty \subset \partial\hAdS$. 

Another example of a hyperbolic plane is the plane dual to an involution itself.
For example,  let $g\in\AdS^3$ be the involution  $\bJ(0)$ defined in \eqref{eq:HyperbolicGeodesic} and corresponding to $i\in\Ht$.
Then one lift of $g^*$ to $\hAdS$ is the unit disc $D^2$ 
in the $xy$-plane defined by $z = 0$. 
The other lift is the exterior of $D^2$ in the $xy$-plane, which is the image
$I_\Ss(D^2)$.
The closure of each disc intersects $\Ss$ in the unit circle in the $xy$-plane. 

An example of a degenerate totally geodesic surface is the Borel subgroup
$\bB = \exp(\bb) \subset \PSLtR$ consisting of the images of 
upper-triangular matrices in $\SLtR$; 
compare the discussion in 
\S\ref{sec:NullPlanes}. Lift $\bB$ to the connected subgroup 
$\hbB\subset\SLtR$, which embeds in $\Ein^3$ by:
\begin{align}\label{eq:UpperTriangulars}
\hbB &\stackrel{\Psi}\longrightarrow \Ein^3 \notag \\
\bmatrix a & b \\ 0 & a^{-1}\endbmatrix
&\longmapsto 
 [a - a^{-1} : b : b : (a + a^{-1}) -2  :  (a + a^{-1}) + 2  ] .
\end{align}
The image is the subspace defined by $Y=Z$ but $U\neq V$.
In affine coordinates $(x,y,z)$ on $\E$, this is the plane $y=z$.
Its closure in $\Ein^3$ is the lightcone $L(p_1)$ where 
\[ p_1 = [0:1:1:0:0] \in C_\infty.\]
Furthermore the boundary of $\Psi(\hbB)$ consists of two rulings
$\Ss_0^+$ (given by $x = 1$) and $\Ss_\pi^-$ 
(given by $x = -1$), as defined in \eqref{eq:Rulings}. 

Another degenerate plane is defined by the subgroup $\bB^\dag$
consisting of the images of {\em lower triangular matrices.\/}
It embeds via
\begin{align}\label{eq:LowerTriangulars}
\hbB^\dag &\stackrel{\Psi}\longrightarrow \Ein^3 \notag \\
\bmatrix a & 0 \\ c & a^{-1}\endbmatrix
&\longmapsto 
 [a - a^{-1} : c : -c : (a + a^{-1}) -2  :  (a + a^{-1}) + 2  ] .
\end{align}
The image is the subspace defined by $Y=-Z$ but $U\neq V$.
In affine coordinates $(x,y,z)$ on $\E$, this is the plane 
$y=-z$.
Its closure in $\Ein^3$ is the lightcone $L(p_2)$ where 
\[ p_2 = [0:1:-1:0:0] \in C_\infty.\]

Furthermore the boundary of $\Psi(\hbB)$ consists of two rulings
$\Ss_0^-$ (given by $x = 1$) and $\Ss_0^-$ 
(given by $x = -1$), as defined in \eqref{eq:Rulings}. 

Finally we consider an indefinite plane discussed in \S\ref{sec:Indefinites}.
It maps to $\Ein^3$ by:
\begin{align}\label{eq:IndefinitePlane}
\mathfrak{S} &\stackrel{\Psi}\longrightarrow \Ein^3 \notag \\
\bmatrix a & b \\ c & a\endbmatrix
&\longmapsto 
 [0: b+ c : b -  c :  2 a  -2  :  2 a + 2  ] 
\end{align}
where $a^2 - b c = 1$.  The image $\Psi(\mathfrak{S})$ is defined
by $X = 0$ and $ U \neq V$.  
%


\section{Crooked surfaces adapted to an involution}\label{sec:AdaptedCS}

This section develops the theory of crooked surfaces, especially
as they relate to an involution $I_\Ss$ in an Einstein hypersphere.
The main goal is to show that  $\AdS$-crooked planes correspond to
crooked surfaces in $\Ein^3$ which are adapted to $I_\Ss$.

A useful invariant of a crooked plane is its {\em spine:\/} the unique spacelike geodesic $\sigma$ lying on the crooked
plane. The spine contains the vertex $p$ and the pair $(p,\sigma)$
determines the crooked plane.

In a similar way, we define a {\em spine\/} $\bars$ for a crooked surface $\CS$.
When $\CS$ is the closure of a crooked plane $\CPM\in\E$,
then $\bars$ equals the union $\sigma \cup \{p_\infty\}$, where 
$p_\infty$ is the improper point. We relate this invariant to the 
parametrization of crooked surfaces by stem-configurations 
discussed in~\cite{BCDGM}.
The spine $\bars$ is a spacelike circle containing the vertex $p_0$
and the improper point $p_\infty$. In general triples $(\bars,p_0,p_\infty)$,
where $\bars$ is a spacelike circle in $\Ein^3$ and $p_0,p_\infty\subset \bars$
are distinct points, provide an alternate parametrization of crooked surfaces.

We say that $\CS$ is {\em adapted\/} to $I_\Ss$ if and only if its spine
$\bars$ is invariant under $I_\Ss$ and $p_0, p_\infty$ are interchanged
under $I_\Ss$.

\subsection{Crooked surfaces.}

A {\em crooked surface\/} is a piecewise smooth hypersurface in $\Ein^3$ 
constructed as follows. 
Take two points  $q_0, q_\infty$ in $\Ein^3$ which do not lie on a photon.
Call one of them, say $q_0$, the {\em vertex\/}.
Call the other one, $q_\infty$, the {\em covertex.\/}
The intersection of the lightcones $L(q_0) \cap L(q_\infty)$
is a spacelike circle in $\Ein^3$ disjoint from $\{q_0,q_\infty\}$.
Choose two distinct points $q_1, q_2 \in\Ein^3$ which lie on 
$L(q_0) \cap L(q_\infty)$. Call the points $q_1,q_2$ the 
{\em hingepoints\/} of the crooked plane.

For $i=0,\infty$ and $j = 1,2$, let $\phi_i^j$ be the unique photon containing
$q_i$ and $q_j$. Call the photons $\phi_0^1, \phi_0^2$ 
{\em hinges\/} and the photons $\phi_\infty^1, \phi_\infty^2$ 
{\em cohinges.\/} 
Call such a quadruple $(q_0,q_\infty;q_1,q_2)$ 
a {\em stem configuration.\/} 
Such a configuration lies on a unique Einstein hypersphere $W$. 

The hinges and cohinges are rulings of $W$ and decompose $W$ into
four squares, each of which have the same four vertices 
$q_0,q_1,q_\infty,q_2$. Two of these squares contain timelike directions emanating 
from $q_0$ or $q_\infty$, and constitute the two components 
$T_1,T_2$ of the {\em stem.\/}

The {\em wings\/} are obtained from the lightcones $L(q_j)$ as follows.
The photons $\phi^j_0,\phi^j_\infty$ separate $L(q_j)$ into two planar bigons, one of which attaches to $T_1\cup T_2$.
For $j=1$, denote these planar bigons by
\[
W_1, W_1' \subset L(q_1j) \setminus (\phi^j_0\cup\phi^j_\infty).\] 
Then exactly one of the planar bigons in $L(q_2)$ 
will be disjoint from $W_1$; call that one $W_2$ and
the other one $W_2'$. Then $W_2'$ will be disjoint from
$W_1'$. The pairs $(W_1,W_2)$ and $W_1', W_2'$
will be the wings of the two crooked surfaces
defined by $(p_0,p_\infty,q_1,q_2)$. 

Define {\em crooked surfaces\/}  
$\CS(p_0,p_\infty; p_1,p_2)$
$\CS(p_0,p_\infty; p_1,p_2)'$ as:
\begin{align*}
\CS(p_0,p_\infty; p_1,p_2) &:= 
T_1 \cup T_2 \cup W_1 \cup W_2 \\
\CS(p_0,p_\infty; p_1,p_2)' &:= 
T_1 \cup T_2 \cup W_1' \cup W_2'. \end{align*}

In \cite{BCDGM}, an ambient orientaton on $\E$
is assumed. This orientation on the Minkowski
patch $\Ein^3\setminus L(p_\infty)$ enables a preference
for ``positively extended'' versus ``negatively extended''
wings (and crooked planes). Since $\Ein^3$ itself is
nonorientable, there is no preferred orientation on $\E$.
Indeed the involution corresponding to the antipodal
map 
\begin{align*} 
\E &\xrightarrow{\mathcal{A}} \E \\
(x,y,z)& \longmapsto (-x,-y,-z) \end{align*}
(given by the orthogonal matrix $-\one_3 \oplus \one_2$)
fixes $p_0,p_\infty,p_1,p_2$ but interchanges the two
crooked surfaces above.
As $\mathcal{A}$ reverses orientation on the Minkowski
patch $\E$, it interchanges the
positively extended crooked planes and 
the negatively extended crooked planes.

\subsubsection{The spine of a crooked surface.\/}

Spacelike circles in $\Ein^3$ correspond to projectivized nullcones of
three-dimensional subspaces $ F \subset \R^{3,2}$ of signature $(2,1)$.
That is, $F$ is the image of an isometric embedding $\R^{2,1} \hookrightarrow
\R^{3,2}$ of inner product spaces. 
Its orthogonal complement $F^\perp$ is nondegenerate as well,
and has signature $(1,1)$ and together $F, F^\perp$ define an orthogonal
decomposition
$$
\R^{2,1} \oplus \R^{1,1} \ \xrightarrow{\cong} 
F \oplus F^\perp \ = \    \R^{3,2}.
$$
The projectivized null cone of $\R^{1,1}$ defines an unordered pair of
non-incindent points in $\Ein^3$.
Orthogonal complementation defines a natural bijection between spacelike
circles and unordered pairs of nonincident points.
When applied to such a pair  $\{p_1,p_2\}\subset\Ein^3$, 
the corresponding spacelike circle equals the intersection of lightcones $L(p_1) \cap L(p_2)$.

Suppose that $\CS = \CS(p_0,p_\infty;p_1,p_2)$ is a crooked surface.
Define the {\em spine\/} $\bars$ of $\CS$ as the spacelike circle dual
to the pair of hingepoints $p_1,p_2$:
$$
\bars\  :=  \  L(p_1) \cap L(p_2).
$$
Since $p_0$ and $p_\infty$ are each incident to $p_i$ for $i=1,2$,
both $p_0$ and $p_\infty$ lie in $\bars$. The triple $(p_0,p_\infty;\bars)$
is equivalent to the stem-configuration $(p_0,p_\infty;p_1,p_2)$ and
equally well determines the crooked surface $\CS$.


\subsubsection{Invariant crooked surfaces and adapted crooked surfaces.\/}
Suppose that $\CS$ is a crooked surface invariant under the involution
$I_\Ss$. Then $I_\Ss$ preserves the stem-configuration and 
permutes the pairs
$p_0, p_\infty$ and $p_1, p_2$. The points $p_0,p_\infty$ are characterized
by the property that the tangent cones of $\CS$ at these points are crooked 
planes. We say that $\CS$ is {\em adapted\/} to $I_\Ss$ if
$I_\Ss$ interchanges $p_0$ and $p_\infty$ but fixes $p_1$ and $p_2$.

This condition is stronger than being invariant. The only spacelike circles
$C$ invariant under $I_\Ss$ meet $\Ss$ in two points. 
Thus the restriction of $I_\Ss$ is an orientation-reversing involution of $C$
fixing this pair. 
In particular $I_\Ss$-invariant crooked surfaces exist where $I_\Ss$ do not
fix the hingepoints, and these will {\em not\/} be adapted in our definition.

\subsection{Completing $\hAdS$ crooked planes to crooked surfaces\/}
Next we show (\eqref{maintheorem2} of the Main Theorem) 
that an $\AdS$-crooked plane lifts to an object in
$\hAdS$ whose closure in $\Ein^3$ is a crooked surface, 
as defined in \cite{Frances1,BCDGM,LareauDussault}. 
Call the  special point corresponding to the improper point
in a crooked surface arising from a crooked plane 
in Minkowski space the {\em covertex.\/}
As shown in \cite{BCDGM}, 
the vertex and covertex are related by 
an automorphism of $\Ein^3$ which preserves the crooked surface.
The condition that the crooked surface
in $\Ein^3$ arises from an $\AdS$-crooked plane is that the vertex and 
covertex are related by the involution $I_\Ss$.

Begin with an $\AdS$-crooked plane $\CPa(g,\bs)$ as in \S\ref{sec:AdSCP}.
Working in the $\PSLtR$-model, we may assume that the vertex $g$ equals
the identity element $e\in\PSLtR$.

Now lift the crooked plane to an $\hAdS$-crooked plane
$\hCPa(\hp,\bs)$ in $\hAdS$ 
as in \S~\ref{sec:LiftCP}. Then the vertex $\hg$ is the identity matrix $\one\in\SLtR$.
The covertex is $-\one\in\SLtR$, the image of $\one$ under the deck involution. We identify the tangent spaces $T_e(\AdS^3)$ with 
$T_{\one}(\hAdS)$
and denote the unit spacelike vector in $T_{\one}(\hAdS)$ corresponding to
$\bs$ by $\hat\bs$.

Under the embedding $\Psi$ of $\hAdS$ in $\Ein^3$,
the vertex $\one$ corresponds to the origin $0 = (0,0,0)$ in Minkowski space $\E$
and $\hat\bs$ corresponds to a unit spacelike vector (still denoted $\hat\bs$)
in $T_0\E \cong \V$. 
By applying yet another isometry of $\AdS^3$, we can assume that $\hat\bs$
is dual to the geodesic 
$\Ho\subset\Ht$ defined by \eqref{eq:HyperbolicGeodesic},
corresponding to the matrices in $\SLtR$ with zero diagonal entries.

Now we build up the picture of the $\hAdS$-crooked plane and show
that its closure in $\Ein^3$ is the crooked surface which is the closure
of a standard crooked plane in $\E$. The main difference are the 
parametrizations of particles and the spine, since these non-null
geodesics are complete in both the realization in Minkowski space
and anti-de Sitter space.  In Minkowski space these geodesics limit
to the improper point $p_\infty$ while in $\AdS^3$ these geodesics limit
to $\Ein^2 = \Ss \cup \mathfrak{I}$.

\subsubsection{Particles and the stem}
We begin with a particle, that is a timelike geodesic, and show that
it is just a reparametrized particle in Minkowski space compactified by the improper point
$p_\infty$.
The stem is a union of particles, whose intersection points with $p^*$
form a geodesic in the  hyperbolic plane $p^*$. 
As in \S\ref{sec:HtInEin}, the two lifts of $p^*$ correspond to the
components of the $2$-sheeted hyperboloid $x^2 + y^2 - z^2 = -1$. 
The geodesic $\sigma\subset p^*$ corresponding to $\Ho\subset\Ht$ is represented by the hyperbola
$$ 
x = y^2 - z^2 + 1 = 0.
$$

The  stem was defined in \S\ref{sec:AdSCP} 
as the union of the timelike geodesics joining
$p$ to $\sigma\subset p^*$. 
Parametrizing the geodesic $\sigma\subset\Ht$
by $\bJ(t)$ (where $t\in\R$), 
as in \eqref{eq:HyperbolicGeodesic}, 
the corresponding elliptic one-parameter subgroups are:
$$
\mathcal{E}(\theta,t) :=
\cos(\theta) \one + \sin(\theta) \bJ(t) = 
\bmatrix   \cos(\theta) & -\sin(\theta) e^t \\ \sin(\theta) e^{-t} & \cos(\theta) \endbmatrix
$$
we obtain the stem of an $\AdS$-crooked plane $\CPa$.
The corresponding subset of $\SLtR$ is defined by $a = d$ with 
$-2 < a = d < 2$. 
Compare  \eqref{eq:StemPlaneInAdS} and \S\ref{sec:Indefinites}.

Under the embedding $\hAdS\hookrightarrow \Ein^3$ 
in \eqref{eq:EmbedSL2R}, 
this is an open subset in the hypersurface $X=0$, that is the
subset defined by:
$$
x = 0, \quad y^2 - z^2 < 0,
$$
which defines the stem of the standard example of a crooked plane, 
as in \cite{DrummGoldman1}. 
The embedding of this totally geodesic submanifold of $\hAdS$ is
explicitly given in \eqref{eq:IndefinitePlane}.

\subsubsection{Hinges and cohinges}\label{sec:HingeCohinge}
The stem is bounded by two hinges, which in 
$\AdS^3$ is the pair of unipotent one-parameter
subgroups
\begin{align}\label{eq:unipotents}
U_1 &:= \pm \bigg\{ \bmatrix 1 & t \\ 0 & 1 \endbmatrix\  \Big|\  t\in\R\bigg\} \notag\\
U_2 &:= \pm \bigg\{ \bmatrix 1 & 0 \\ t & 1 \endbmatrix\  \Big|\  t\in\R\bigg\} 
\end{align}
which respectively correspond to the {\em hinges\/} $h_1, h_2\subset\AdS^3$.
By \eqref{eq:unipotents}, their $\hAdS$-lifts map via $\Psi$ to arcs of photons
\begin{align}\label{eq:hinges}
\hat{h}_1^+ &:= 
\bigg\{ [0:t:t:0:4] \ \Big|\  t\in\R \bigg\} 
 \notag \\
\hat{h}_1^- &:= 
\bigg\{  [0:t:t:4:0] \ \Big|\  t\in\R \bigg\} \notag \\
\hat{h}_2^+ &:= 
\bigg\{ [0:t:-t:0: 4] \ \Big|\  t\in\R \bigg\} \notag \\
\hat{h}_2^- &:= 
\bigg\{ \ [0:t:-t:4:0] \ \Big|\  t\in\R \bigg\} 
\end{align}
in homogeneous coordinates, respectively.
The hinges $\hat{h}_i$ of the tangent cone $\CPM$
exponentiate to the hinges of the $\hAdS$ crooked plane
$\hCPa$. Since the exponential map is linear (with speed
$1/4$), the hinges of $\hCPa$ are also complete geodesics
in $\hAdS$. Their closures consist of their unions
with hingepoints which lie on $C_\infty$.
The cohinges $\hat{h}_i^-$ are just the images of the hinges
under inversion $I_\Ss$, and their closures also
consist of their unions with hingepoints:

\begin{align*}
\partial\hat{h}_1^+ \ =\  \partial\hat{h}_1^- \ = p_1 
& =  \  [0:1:1:0:0] \\
\partial\hat{h}_2^+  \ = \ \partial\hat{h}_2^- \ = p_2
& =  \ [0:1:-1:0:0],
\end{align*}
respectively. 
In Minkowski space these hinges are just the lines
consisting of points $(0,y,z)\in \E$ where
$y = z$ and $y= -z$ respectively,
and the cohinges are their $I_\Ss$-images in the inverse cone $\mathfrak{I}$.
The hingepoints  $p_i = \partial\hat{h}_1^{\pm}$ (for $i=1,2$)
each lie on the ideal circle $C_\infty$.

\subsubsection{Wings and cowings}
From the hinges $h_1,h_2$ we determine the halfplanes 
$h_i^{\perp}$  (where $i=1,2$) in Minkowski space. 
These will intersect $\hAdS$  in totally geodesic null surfaces, 
and the complements of the hinges in these planes will be halfplanes, 
one of which will be the wings of the crooked plane.
As the hinges $h_\pm$ correspond to unipotent one-parameter subgroups 
$U_\pm$ (as in \eqref{eq:unipotents} in \S\ref{sec:HingeCohinge}), 
the associated null planes $h_\pm^{\perp}$ correspond
to the normalizers of these subgroups $B_\pm := \Norm(U_\pm)$:
\begin{equation}\label{eq:BorelSubgroups}
\pm \bmatrix a & b \\ 0 & a^{-1}\endbmatrix,
\pm \bmatrix a & 0 \\ c & a^{-1} \endbmatrix
\end{equation}
which define totally geodesic surfaces in $\hAdS$ which are degenerate.
The complement $h_{\pm}^{\perp} \setminus h_\pm$ 
is a disjoint union of two null halfplanes.

The unipotent one-parameter subgroup $U$ 
(which is the commutator subgroup of $B_\pm$)
is defined by $a = 1$. While $U$ is the stabilizer of the vector 
$$
\bmatrix 1 \\ 0 \endbmatrix \in \R^2,
$$
its normalizer $B$ is the stabilizer of the corresponding point 
(denoted $\infty$)  in 
the projective line $\P(\R^2)$.
Furthermore $U$ separates $B$ into two subsets depending on
whether $a> 1$ or $a < 1$. The complement $B\setminus U$ consists
of hyperbolic elements in $B$, and the two subsets are distinguished by
whether $\infty$ is an attracting (respectively repelling) fixed point.

The corresponding subsets  in $\Ein^3$ are, respectively:
$$ [a - a^{-1} : b : b : (a + a^{-1}) -2  :  (a + a^{-1}) + 2  ] $$
$$ [a - a^{-1} : c : -c : (a + a^{-1}) -2  :  (a + a^{-1}) + 2  ] $$
which are defined in homogeneous coordinates by $Y = Z $ and
$ Y = -Z$ respectively, and in affine Minkowski coordinates
by $y=z$ and $y=-z$ respectively.

Consider the first subset $Y=Z$ which contains the hinge $X = U =  Y - Z = 0$.
The complement of the hinge in the wing is a pair of null half-planes, 
each of which intersects $\AdS^3$ in a strip bounded by the hinge and a null line
which is a ruling of $\Ss$.

\subsubsection{The spine}
As the wings in our examples correspond to subsemigroups of 
Borel subgroups consisting of upper- and lower-triangular matrices,
respectively, the spine consists of the Cartan subgroup of diagonal matrices.
In particular \eqref{eq:Hyperbolics} in \S\ref{sec:Geodesics} describes
the mapping from the spine to the interval \eqref{eq:UnitInterval} on the
$x$-axis, 
and \eqref{eq:MinusHyperbolics} describes
the mapping from the cospine to the union of two rays 
\eqref{eq:TwoRays} on the $x$-axis, separated by the improper point at
infinity.
The spine $\bars$ of the corresponding crooked plane is then the $x$-axis
completed by the improper point.

\subsection{Characterization of $\AdS$-crooked planes\/}
Finally we prove the converse statement \eqref{maintheorem3} of the Main Theorem.

Let $(q_0,q_\infty;q_1,q_2)$ be a stem configuration and 
$$
\CS \ := \ \CS(q_0,q_\infty;q_1,q_2) \ \subset \ \Ein^3
$$
the corresponding crooked surface.
Suppose that $\CS$ is $I_\Ss$-adapted.
Then $I_\Ss$ interchanges the vertex $q_0$ and the covertex $q_\infty$
and fixes the hingepoints $q_1$ and $q_2$.

Then,  since $I_\Ss$ fixes neither $q_0$ nor $q_\infty$,
both points lie in $\Psi(\hAdS)$. Since the group of isometries
$\mathsf{Isom}(\hAdS)$ acts transitively
on $\hAdS$,  we may apply 
an automorphism of $\AdS^3$ to assume that
the covertex of $\CPa$ is the improper point, that is,
$q_\infty = p_\infty$.  Since $I_\Ss$ interchanges the
origin and the improper point, and 
interchanges the covertex and the vertex 
(because $\CPa$ is $I_\Ss$-adapted),
the vertex equals the origin: 
$q_0 = p_0$.  
Because the covertex of $\CS$ equals $p_\infty$,
the crooked surface is the closure of a crooked plane
$\CPM$ in Minkowski space:
$$
\CPM\ = \ \CS \setminus L(p_\infty) \ = \ \CS \cap \E
$$
and the vertex of this crooked plane is the origin $p_0$.

Now we adjust the hingepoints.
The group of orientation-preserving 
isometries $\hAdS$ which fix $\one$ and $-\one$ 
is just the group $H$ of inner automorphisms of $\SLtR$, 
which is isomorphic to the identity component of $\SOto$.
Similarly the group of orientation-preserving isometries of
$\E$ which fix the origin equals $\SOto$, and the embedding
$\Psi$ is equivariant with respect to the isomorphism
$$
\Inn(\SLtR) \xrightarrow{\cong} \SOto^0 \cong \Isom(\Ht)^0.
$$
Since $\Isom(\Ht)^0$ acts transitively on geodesics in $\Ht$, 
the group $\SOto$ acts transitively on ordered pairs of distinct
null lines in $\Rto$. Therefore there exists an isometry of $\hAdS$
which maps the hingepoints $q_i$ to $p_i$ for $i=1,2$ respectively.
Apply this isometry to assume that $q_i = p_i$. 
Thus $\CS$ is the closure of the double cover $\hCPa$ of an $\AdS$-crooked plane  $\CPa$.
The proof is complete.
\section*{Appendix:\ Einstein hyperspheres bounding $\AdS^3$}
Because $\PSLtR$ is homeomorphic to an open solid torus 
$D^2\times S^1$, its boundary is a
$2$-torus $S^1\times S^1$. 

Since this theory seems not to be so well-known, we begin with a
brief exposition of {\em singular projective transformations.\/}
This is a simple case of the theory of {\em wonderful compactifications\/}
developed by de Concini-Procesi~\cite{DeConciniProcesi}.

A natural compactification of  the group $\PGLtR$ is the projective space
$\rpthree$ obtained by projectivizing the embedding
\[ \GLtR \hookrightarrow \Mat_2(\R)\setminus\{ \zero\} \]
of the group of {\em invertible\/} $2\times 2$ real matrices in the set of
{\em nonzero\/} $2\times 2$ real matrices.
The induced map 
\[ \PGLtR \hookrightarrow \rpthree \]
embeds $\PGLtR$ as an open dense subset of $\rpthree$.
The complement $\rpthree\setminus\PGLtR$ consists
of {\em singular projective transformations\/} of $\rpo$.
It naturally identifies with $\rpo\times\rpo$ and is homeomorphic to a torus.

A point in the the complement $\rpthree\setminus\PGLtR$ is the projective
equivalence class of a $2\times 2$ real matrix $M$ 
which is nonzero but singular.
In other words,  $M$ is a $2\times 2$ matrix of rank $1$. 
The corresponding linear
endomorphism of $\R^2$ has a $1$-dimensional kernel and a $1$-dimensional
range, which may coincide. This pair of lines in $\R^2$ determines
a point
\[ \big(\Ker(M), \Image(M)\big)\ \in \ \rpo\times\rpo. \]
Two matrices giving the same pair are projectively equivalent.
The corresponding {\em singular projective transformation\/} 
is only defined on $\rpo\setminus \{\Ker(M)\}$ where it is the constant mapping 
\[ \rpo\setminus \{\Ker(M)\} \longrightarrow \{\Image(M)\}. \]

Here are two basic examples.
The one-parameter semigroup of projective transformations
defined by the diagonal matrix 
\[ \bmatrix e^t & 0 \\ 0 & e^{-t} \endbmatrix \]
converges,
as $t\longrightarrow + \infty$, 
to the  singular projective transformation defined by the diagonal matrix
\[ \bmatrix 1 & 0 \\ 0 & 0 \endbmatrix \]
which corresponds to the ordered pair 
\[ 
\big([0:1],[1:0] \big) \ \in \ \rpo\times\rpo.\]

Likewise the one-parameter semigroup of projective transformations
\[ \bmatrix 1 & t  \\ 0 & 1 \endbmatrix \]  converges,
as $t\longrightarrow + \infty$,  to the 
singular projective transformation defined by 
\[ \bmatrix 0 & 1 \\ 0 & 0 \endbmatrix, \]
which corresponds to the ordered pair
$$ \big([1:0],[1:0]) \big) \ \in \ \rpo\times\rpo. $$
For more details see \cite{PGOM} and the references cited there.

The group $\PSLtR$ is one component of $\PGLtR$; 
the other component corresponds
to matrices with negative determinant. 
Denote the closure of
$\PSLtR \subset \rpthree$ by $\overline{\PSLtR}$.
Since the inclusion
\[ \PSLtR \hookrightarrow\overline{\PSLtR} \]
is a homotopy-equivalence, 
the compactification of $\PSLtR$ lifts to a compactification of the double covering 
\[ \SLtR \longrightarrow \PSLtR. \] 

Points on the boundary $\partial\SLtR \approx S^1\times S^1$
can be geometrically interpreted as follows. 
Choose an orientation $O_{\R^2}$ on $\R^2$. 
The double covering $\widehat{\rpo}$ of $\rpo$ consists of 
{\em oriented lines,\/} that is, pairs $(\ell,O_\ell)$,
where $\ell\in\rpo$ is a line in $\R^2$ and $O_\ell$ is an orientation
on $\ell$. 
Then the boundary $\partial\SLtR$ consists of ordered pairs of oriented lines
\[ \big( (\ell_1,O_1),\  (\ell_2,O_2)\big) \ \in \ 
\widehat{\rpo}\times \widehat{\rpo} \]
such that:
\begin{itemize}
\item $O_1 = O_2$ if $\ell_1 = \ell_2$;
\item $
O_{\R^2} = O_1 \oplus O_2 $ with respect to the decomposition
$\R^2 = \ell_1 \oplus \ell_2$ if  $\ell_1 \neq \ell_2$.
\end{itemize}
Alternatively, the double covering 
\[ \partial\SLtR \longrightarrow \partial\PSLtR \approx
\rpo \times \rpo\]
is defined by the homomorphism
$\pi_1(\rpo \times \rpo) \longrightarrow \Z/2$ which is nontrivial on each Cartesian factor.

These two projections correspond to the {\em rulings\/}
of $\partial\hAdS$, as described in \S\ref{sec:Rulings}.

 \bibliographystyle{amsplain}

\begin{thebibliography}{99}
%
\bibitem{BBDGGKKSZ}
Barbot, T., Bonsante, F.,  Danciger,  J., 
Goldman, W., Gu\'eritaud, F., Kassel, F.,
Krasnov, K., Schlenker, J.-M.,  Zeghib, A.,
{\em Some open questions on anti-de Sitter geometry. \/}
{\tt arXiv:1205.6103 }

\bibitem{BCDGM}
Barbot, T., Charette, V., Drumm, T., Goldman, W., and Melnick, K.,
{\em A Primer on the (2+1)-Einstein Universe,\/}
in {\em Recent Developments in Pseudo-Riemannian Geometry,\/} 
(D.\ Alekseevsky, H. Baum, eds.) Erwin Schr\"odinger Lectures in Mathematics and
Physics, European Mathematical Society (2008), 179--230
{\tt math.DG.0706.3055}.

\bibitem{BCDG}
Burelle, J.-P., Charette, V., Drumm, T., and Goldman, W.,
{\em Crooked halfspaces,\/}
Enseignement  Math\'ematique (to appear) {\tt math.DG.1211.4177}

\bibitem{CFL}
Charette, V., Francoeur, D., Lareau-Dussault, R.,
{\em Fundamental polyhedra in the Einstein universe,\/}
Topology and its Applications  {\bf 174} (2014), 62--80.
{\tt math.DG.1307.6531}

\bibitem{Danciger}
Danciger, J., 
{\em Geometric transitions from hyperbolic to $\AdS$ geometry,\/}
Doctoral thesis, Stanford University (2011)
\bibitem{DGK}
Danciger, J., Gu\'eritaud, F., and Kassel, F.
{\em Fundamental domains for free groups acting on anti-de Sitter $3$-space,\/} (in preparation)

\bibitem{DeConciniProcesi}
De Concini, C.; Procesi, C.,
{\em Complete symmetric varieties,\/} in Gherardelli, Francesco,
``Invariant theory (Montecatini, 1982)'', 
Lecture Notes in Mathematics, 996, Berlin, New York:  
Springer-Verlag, pp. 1Ð44. (1983)

\bibitem{Drumm}
Drumm, T.,
{\em Fundamental polyhedra for Margulis space-times,\/}
Topology {\bf 31} (4) (1992), 677-683.

\bibitem{DrummGoldman1}
    Drumm, T. and Goldman, W.,
    {\em The geometry of crooked planes,\/}
    Topology {\bf 38}, No. 2, (1999) 323--351.

\bibitem{Eschenburg}
Eschenburg, J.- H.,
{\em Lecture Notes on Symmetric Spaces,\/} 
available from 
{\tt www.math.uni-augsburg.de/eschenbu/symspace.pdf}.

\bibitem{Frances1} Frances, C., 
{\em G\'eometrie et dynamique lorentziennes conformes,\/}
Th\`ese, E.N.S. Lyon (2002)
\bibitem{Frances2}
Frances, C.,
{\em The conformal boundary of Margulis space-times,\/}
C. R. Math. Acad. Sci. Paris 336 (2003), no. 9, 751--756.

\bibitem{PGOM}
Goldman, W.,
{\em Projective Geometry on Manifolds,\/}
Notes from a graduate course at the University of Maryland (1988),
available from 
{\tt http://www.math.umd.edu/$\sim$wmg/pgom.pdf}

\bibitem{Goldman}
Goldman, W.,
``Complex Hyperbolic Geometry,''
Oxford University Press (1999)

\bibitem{Helgason}
Helgason, S.,
``Differential Geometry, Lie Groups, and Symmetric Spaces,''
Pure and Applied Mathematics, a series of monographs and textbooks,
{\bf 80} Academic Press, New York (1978)

\bibitem{LareauDussault}
Lareau-Dussault, R.,
{\em Les surfaces croches de l'universe d'Einstein,\/}
m\'emoire, Universit\'e de Sherbrooke (2012)

\bibitem{Mess}
Mess, G.,
{\em Lorentz spacetimes of constant curvature,\/}
Geometriae Dedicata  {\bf126} (2007), 3Ð45.

\bibitem{O'Neill}
O'Neill, B.,
``Semi-Riemannian Geometry, with Applications to Relativity,''
Pure and Applied Mathematics {\bf 103}, Academic Press, Inc.,
San Diego, California (1983)

\bibitem{Rossi}
Rossi Salvemini, C.,
{\em Espace-temps globalement hyperboliques conform\'ement plats,\/}
Th\`ese, Universit'e de Avignon (2012)

\bibitem{Wolf}
Wolf, Joseph A.,
``Spaces of constant curvature.''
Sixth edition. AMS Chelsea Publishing, Providence, RI, 2011. xviii+424 

\end{thebibliography}

\end{document}